\documentclass[reqno,b5paper]{amsart}
\usepackage{amsmath}
\usepackage{amssymb}
\usepackage{amsthm,amsxtra}
\usepackage{cite}
\usepackage{enumerate}
\usepackage[mathscr]{eucal}
\usepackage{float}
\newtheorem{Def}{Definition}[section]
\newtheorem{Th}{Theorem}[section]
\newtheorem{Ex}{Example}[section]

\newtheorem{Lemma}{Lemma}[section]
\newtheorem{Prop}{Proposition}[section]
\newtheorem{Cor}{Corollary}[section]
\newtheorem{Rem}{Remark}[section]

{}

\newcommand{\kc}{\mathcal{K}}
\newcommand{\rb}{\mathbb{R}}
\newcommand{\nb}{\mathbb{N}}

\newcommand{\rc}{\mathcal{R}}

\newcommand{\ac}{\mathcal{A}}
\newcommand{\bc}{\mathcal{B}}
\newcommand{\bk}{\mathfrak{B}}
\newcommand{\bb}{\mathfrak{b}}
\newcommand{\df}{\mathfrak{d}}

\newcommand{\fc}{\mathcal{F}}
\newcommand{\mc}{\mathcal{M}}
\newcommand{\oc}{\mathcal{O}}
\newcommand{\pc}{\mathcal{P}}
\newcommand{\pf}{\mathfrak{p}}
\newcommand{\qc}{\mathcal{Q}}

\newcommand{\tc}{\mathcal{T}}
\newcommand{\uc}{\mathcal{U}}

\newcommand{\vc}{\mathcal{V}}
\newcommand{\wc}{\mathcal{W}}

\newcommand{\zc}{\mathcal{Z}}

\newcommand{\Ga}{\Gamma}

\newcommand{\obs}{\oc_{\bk^s}}
\newcommand{\bs}{{\bk^s}}
\newcommand{\fbs}{{f(\bk)^s}}
\newcommand{\gbs}{\gamma_{\bk^s}}
\newcommand{\Gbs}{\Gamma_{\bk^s}}
\newcommand{\obhs}{\oc_{(\widehat{\bk})^s}}
\newcommand{\Gbhs}{\Gamma_{(\widehat{\bk})^s}}

\DeclareMathOperator{\Sf}{{\sf S}_{fin}}
\DeclareMathOperator{\S1}{{\sf S}_1}
\DeclareMathOperator{\Uf}{{\sf U}_{fin}}
\DeclareMathOperator{\Gf}{{\sf G}_{fin}}
\DeclareMathOperator{\G1}{{\sf G}_1}
\DeclareMathOperator{\cov}{{\sf cov}}
\DeclareMathOperator{\add}{{\sf add}}
\newcommand{\bp}{\begin{proof}}
\newcommand{\ep}{\end{proof}}

\newcommand{\obns}{\oc_{(\bk^n)^s}}
\newcommand{\bns}{(\bk^n)^s}
\newcommand{\bn}{\bk^n}

\newcommand{\cdr}{{\sf CDR}_{\sf sub}}
\begin{document}

\title[Further observations on bornological covering properties]{Further observations on bornological covering properties and selection principles}

\author[D. Chandra, P. Das  and S. Das]{Debraj Chandra$^{\dag}$, Pratulananda Das$^*$  and Subhankar Das$^*$\ }

\address{\llap{$\dag$\,}Department of Mathematics, University of Gour Banga, Malda-732103, West Bengal, India}
\email{debrajchandra1986@gmail.com}

\address{\llap{*\,}Department of Mathematics, Jadavpur University, Kolkata-700032, West Bengal, India}
\email{pratulananda@yahoo.co.in, subhankarjumh70@gmail.com}

\thanks{The third author
is thankful to University Grants Commission (UGC), New Delhi-110002, India for granting UGC-NET Junior Research Fellowship (1183/(CSIR-UGC NET DEC.2017)) during the tenure of which this work was done.}


\begin{abstract}
This article is a continuation of the study of bornological open covers and related selection principles in metric spaces done in (Chandra et al. 2020 \cite{dcpdsd}) using the idea of strong uniform convergence (Beer and Levi, 2009 \cite{bl}) on bornology. Here we explore further ramifications, presenting characterizations  of various selection principles related to certain classes of bornological covers using the Ramseyan partition relations, interactive results between the cardinalities of bornological bases and certain selection principles involving  bornological covers, producing new observations on the $\bs$-Hurewicz property introduced in \cite{dcpdsd} and several results on the $\bs$-Gerlits-Nagy property of $X$ which is introduced here following the seminal work of \cite{gn}. In addition, in the finite power $X^n$ with the product bornology $\bk^n$, the $\bns$-Hurewicz property as well as the $\bns$-Gerlits-Nagy property of $X^n$ are characterized in terms of properties of $(C(X),\tau^s_\bk)$ like countable fan tightness, countable strong fan tightness along with the Reznichenko's property.
\end{abstract}
\maketitle
\smallskip
\noindent{\bf\keywordsname{}:} {Bornology, selection principles, open $\bs$-cover, $\gamma_\bs$-cover, $\bs$-Hurewicz property, $\bs$-Gerlits-Nagy property, Ramseyan partition relation, topology of strong uniform convergence, function space $C(X)$.}

\section{Introduction}

This article is continuation of our work in \cite{dcpdsd} and like that paper, here also we follow the notations and terminologies of \cite{arh,mcnt,engelking,hh}. The primary structure where we would work is a bornology.  A bornology $\bk$ on a metric space $(X, d)$ is a family of subsets of $X$ that is closed under taking finite unions, is hereditary and forms a cover of $X$ (see \cite{hh}). A base $\bk_0$ for a bornology $\bk$ is a subfamily of $\bk$ that is cofinal in $\bk$ with respect to inclusion i.e.  for $B\in \bk$ there is a $B_0\in \bk_0$ such that $B\subseteq B_0$. A base is called closed (compact) if all of its members are closed (compact) which would be most useful in our endevors. We list a few natural bornologies on $X$ as follows: $(1)$ The family $\fc$ of all finite subsets  of $X$, the smallest bornology on $X$; $(2)$ The family of all non empty subsets of $X$, the largest bornology on $X$; $(3)$ The family of all non empty $d$-bounded subsets of $X$; $(4)$ The family $\kc$ of non empty subsets of $X$ with compact closure (more examples can be seen from \cite{dcpdsd}).

In the bornological investigations, the notion of continuity which turned out to be the most useful is the notion of strong uniform continuity on a bornology (introduced in \cite{bl} by Beer and Levi). A mapping $f : X\rightarrow Y$ where $(Y,\rho)$ is another metric space, is strongly uniformly continuous on a subset $B$ of $X$ if for each $\varepsilon>0$ there is a $\delta >0$ such that $d(x_1, x_2)< \delta$ and $\{x_1, x_2\}\cap B \neq \emptyset$ imply $\rho( f(x_1),f(x_2))<\varepsilon $.
For a bornology $\bk$ on $X$, $f$ is called strongly uniformly continuous on $\bk$ if $f$ is strongly uniformly continuous on $B$ for each $B \in\bk$. They had also introduced a new topology on $Y^X$ the set of all functions from $X$ into $Y$, called the topology of strong uniform convergence and studied various properties in function spaces. This study has been further continued in \cite{cmh} .

In \cite{cooc1} (see also \cite{cooc2}), M. Scheepers began a systematic study of selection principles in topology and their relations to game theory. Those interested in the long and illustrious history of selection principles and its recent developments can consult the papers \cite{sur1,sur2,tb} where many more references can be found. Using the topology of strong uniform convergence on a bornology, the study of open covers and related selection principles in function spaces had been initiated in \cite{cmk}. In \cite{dcpdsd}, we carried out further advancement in this direction where the main focus was to obtain Schepeers' like diagrams and study the notion of strong-$\bk$-Hurewicz property (or $\bs$-Hurewicz property).

This paper intends to complete the line of investigations started in \cite{dcpdsd} and is organised as follows. In section $2$, we present some basic observations on bornological covers and related selection principles under continuous functions. In section $3$,  we first investigate the behaviour of certain selection principles involving bornological covers under the cardinality of a base $\bk_0$ and then obtain their characterizations in terms of Ramseyan partition relations. Section $4$ is the most important part of this article where we first focus more on the $\bs$-Hurewicz property. Considering the product bornology $\bk^n$ on $X^n$, the $\bns$-Hurewicz property of $X^n$ is shown to be equivalent with the $\bs$-Hurewicz property of $X$ and moreover it is characterized Ramsey-theoretically as well as game-theoretically. Later in this section, we introduce $\bs$-Gerlits-Nagy property of $X$ (following the seminal work of \cite{gn}) and proceed to establish equivalence with $\bns$-Gerlits-Nagy property of $X^n$ and further use Ramseyan partition relations to characterize it. Finally in Section $5$, we devote our attention to the function space $C(X)$ equipped with the topology of strong uniform convergence $\tau^s_\bk$ on $\bk$, where our primary objective is to present characterizations of $\bns$-Hurwicz property as well as $\bns$-Gerlits-Nagy property of $X^n$ in terms of properties of $(C(X),\tau^s_\bk)$ like countable fan tightness, countable strong fan tightness along with the Reznichneko's property. In addition, we also present some characterizations of $\alpha_i$ properties using selection principles.

\section{Preliminaries}

We follow the notations and terminologies of \cite{arh, engelking, hh, mcnt}. Throughout the paper $(X,d)$ stands for an infinite metric space and $\nb$ stands for the set of positive integers.
We first write down two classical selection principles formulated in general form in \cite{cooc1,cooc2}.
For two nonempty classes of sets $\ac$ and $\bc$ of an infinite set $S$, we define

 $\S1(\ac,\mathcal B)$: For each sequence $\{A_n:n\in \nb \}$ of elements of $\ac$, there is a sequence $\{b_n:n\in \nb\}$ such that $b_n\in A_n$ for each $n$ and $\{b_n:n\in \nb\}\in \mathcal B$.

$\Sf(\ac,\mathcal B)$: For each sequence $\{A_n:n\in \nb\}$ of elements of $\ac$, there is a sequence $\{B_n:n\in \nb\}$ of finite  (possibly empty) sets such that $B_n\subseteq A_n$ for each $n$ and $\bigcup_{n\in \nb}B_n\in \mathcal B$.

There are infinitely long games corresponding to these selection principles.

$\G1(\ac,\bc)$ denotes the game for two players, ONE and TWO, who play a round for each positive integer $n$. In the $n$-th round ONE chooses a set $A_n$ from $\ac$ and TWO responds by choosing an element $b_n\in  A_n$. TWO wins the play $\{A_1,b_1, \dotsc, A_n, b_n, \dotsc \}$ if $\{b_n :n\in \nb\}\in \bc$. Otherwise ONE wins.

$\Gf(\ac,\bc)$ denotes the game where in the $n$-th round ONE chooses a set $A_n$ from $\ac$ and TWO responds by choosing a finite (possibly empty) set $B_n\subseteq A_n$. TWO wins the play $\{A_1,B_1, \dotsc,A_n,B_n, \dotsc  \}$ if $\bigcup_{n\in \nb}B_n\in \bc$. Otherwise ONE wins.

We will also consider

$\Uf(\ac,\mathcal B)$: For each sequence $\{A_n:n\in \nb\}$ of elements of $\ac$, there is a sequence $\{B_n:n\in \nb\}$ of finite  (possibly empty) sets such that $B_n\subseteq A_n$ for each $n$ and either $\{\cup B_n:n\in \nb\}\in \mathcal B$ or for some $n$, $\cup B_n=X$ (from \cite{cooc1,cooc2}).

$\cdr(\ac,\bc)$: For each sequence $\{A_n:n\in \nb\}$ of elements of $\ac$ there is a sequence $\{B_n:n\in \nb\}$ such that for each $n$, $B_n\subseteq A_n$, for $m\neq n$, $B_m\cap B_n=\emptyset$ and each $B_n$ is a member of $\bc$\cite{cooc1}.

The following selection principles called $\alpha_i$ properties, are defined in \cite{alpha}. The symbol $\alpha_i(\ac,\bc)$ for $i=1,2,3,4$ denotes that for each sequence $\{A_n:n\in \nb\}$ of elements of $\ac$, there is a $B\in \bc$ such that

$\alpha_1(\ac,\bc)$: for each $n\in \nb$, the set $A_n\setminus B$ is finite.

$\alpha_2(\ac,\bc)$: for each $n\in \nb$, the set $A_n\cap B$ is infinite.

$\alpha_3(\ac,\bc)$: for infinitely many $n\in \nb$, the set $A_n\cap B$ is infinite.

$\alpha_4(\ac,\bc)$: for infinitely many $n\in \nb$, the set $A_n\cap B$ is non empty.\\

Now we recall definitions of some symbols related to Ramseyan partition relation.
Let $\uc$ be an element in $\ac$. A function $f:[\uc]^2\rightarrow \{1,\dotsc,k\}$ is said to be a coloring \cite{cooc1} if for each $U\in \uc$ and every $\vc\in \ac$ with $\vc\subseteq\uc$, there is a $i\in \{1,\dotsc,k\}$ such that $\{V\in \vc:f(\{U,V\})=i\}$ belongs to $\ac$.

We say that $X$ satisfies the partition relation $\ac\rightarrow \lceil\bc\rceil^2_k$ for $k\in \nb$ if for every $\uc\in \ac$ and any coloring $f:[\uc]^2\rightarrow \{1,\dots,k\}$ there are a $i\in \{1,\dots,k\}$, a set $\vc\in \bc$ with $\vc\subseteq\uc$ and a finite to one function $\phi:\vc\rightarrow \nb$ such that for every $V,W\in \vc$ with $\phi(V)\neq \phi(W)$, $f(\{V,W\})=i$ \cite{betd, cooc1}. In this case $\vc$ is said to be eventually homogeneous for $f$. This symbol is known as Baumgartner–Taylor partition symbol.

We say that $X$ satisfies the partition relation $\ac\rightarrow (\bc)^n_k$ for $n,k\in \nb$  if for every $\uc\in \ac$ and any coloring $f:[\uc]^n\rightarrow \{1,\dots,k\}$ there are a $i\in \{1,\dots,k\}$ and a set $\vc\in \bc$ with $\vc\subseteq \uc$ such that for each $V\in [\vc]^n$, $f(V)=i$ \cite{Ramsey}. Also in this case $\vc$ is said to be homogeneous for $f$. This symbol is known as ordinary partition symbol.\\

The following cardinal numbers (see \cite{vaughan} for more details) will be used in the sequel.
The eventually dominating order $\leq^*$ on the Baire space $\nb^\nb$ is defined as follows. For  $f,g\in \nb^\nb$, we say that $f\leq^* g$ if $ f(n)\leq g(n)$ for all but finitely many $n$. Let $A$ be a subset of $\nb^\nb$. The set $A$ is bounded if there is a function $ g\in\mathbb{N}^\mathbb{N}$ such that $f\leq^* g$ for all $f\in A$. The symbol $\bb$  denotes the minimal cardinality of an unbounded subset of $(\nb^\nb,\leq^*)$.
A subset $A$ of $\nb^\nb$ is dominating if for each function $g\in \nb^\nb$ there exists a  function $f\in A$ such that $g\leq^* f$. The symbol $\df$  denotes the minimal cardinality of a dominating subset of $(\nb^\nb,\leq^*)$.
Let $\ac$ be a family of infinite subsets of $\nb$. $P(\ac)$ denotes that there is a subset $P$ of $\nb$ such that for each $A\in \ac$, $P\setminus A$ is finite. The symbol $\pf$ denotes the smallest cardinal number $k$ for which the following statement is false:
For each family $\ac$ if any finite subfamily of $\ac$ has infinite intersection and $|\ac|\leq k$ then $P(\ac)$ holds. The symbol $\cov(\mc)$  denotes the smallest cardinal number $k$ such that
a family of $k$ first category subsets of the real line covers the real line. For any family $\ac$ of subsets of $\nb^\nb$ with cardinality less than $\cov(\mc)$ implies that there is a $g\in \nb^\nb$ such that for every $f\in \ac$ the set $\{n\in \nb:f(n)=g(n)\}$ is infinite.
The symbol $\add(\mc)$ denotes the smallest cardinal number $k$ such that there is a
family of $k$ first category sets of real numbers whose union is no longer first category.
The following relations between the cardinal numbers mentioned above
are  well known. $\add(\mc)\leq \cov(\mc)\leq \df$, $\pf\leq \bb\leq \df$, $\pf\leq \cov(\mc)$ and also $\add(\mc)=\min\{\bb,\cov(\mc)\}$.\\

For $x\in X$, we denote $\Omega_x=\{A\subseteq X:x\in \overline{A}\setminus A\}$ \cite{coocvii}.
$X$ is said to have the Reznichenko  property at $x\in X$ if for each countable set $A$ in $\Omega_x$ there is a partition $\{A_n:n\in \nb\}$ of $A$ into pairwise disjoint finite subsets of $A$ such that for each neighbourhood $W$ of $x$, $W\cap A_n\neq \emptyset$ for all but finitely many $n$ \cite{kocrez, kocfs}. The collection of all such countable sets is denoted by $\Omega^{gp}_x$. $X$ is said to be Fr\'{e}chet–Urysohn if for each  subset $A$ of $X$ and each $x\in A$ there is a sequence in $A$ converging to $x$. $X$ is strictly Fr\'{e}chet Urysohn (in short SFU) if $\S1(\Omega_x, \Sigma_x)$ holds for each $x\in X$. $X$ is said to have countable tightness if for every $x\in X$ and $A\in \Omega_x$ there is a countable subset $B$ of $A$ such that $B\in \Omega_x$  \cite{arh}.  Also $X$ is said to have countable fan tightness (countable strong fan tightness)  at $x$ if $X$ satisfies $\Sf(\Omega_x, \Omega_x)$ ($\S1(\Omega_x, \Omega_x)$) \cite{arh,sakai}. The symbol $\Sigma_x$ denotes the collection of all sequences that converges to $x\in X$ \cite{kock2}.\\

Next we recall some classes of bornological covers of $X$. Let $\bk$ be a bornology on the metric space $(X,d)$ with closed base. For $B\in \bk$ and $\delta>0$, let $B^\delta= \bigcup_{x\in B}S(x,\delta)$, where $S(x,\delta)=\{y\in X:d(x,y)<\delta\}$. It can be easily checked that $\overline{B^\delta}\subseteq B^{2\delta}$ for every $B\in \bk$ and $\delta>0$. A cover $\uc$ is said to be a strong-$\bk$-cover (in short, $\bs$-cover)\cite{cmh} if $X\not\in \uc$ and for each $B\in \bk$ there exist $U\in \uc$ and $\delta>0$ such that $B^\delta\subseteq U$. If the members of $\uc$ are open then $\uc$ is called an open $\bs$-cover. The collection of all open $\bs$-covers is denoted by $\oc_\bs$. An open cover $\uc=\{U_n:n\in \nb\}$ is said to be a $\gbs$-cover \cite{cmh} (see also \cite{cmk}) of $X$, if it is infinite and for every $B\in \bk$ there exist a $n_0\in \nb$ and a sequence $\{\delta_n:n\geq n_0\}$ of positive real numbers satisfying $B^{\delta_n}\subseteq U_n$ for all $n\geq n_0$. The collection of all $\gbs$-covers is denoted by $\Gbs$. An open cover $\uc$ of $X$ is said to be $\bk^s$-groupable \cite{dcpdsd} if it can be expressed as a union of countably many finite pairwise disjoint sets $\uc_n$ such that for each $B\in \bk$ there exist a $n_0\in \nb$ and a sequence $\{\delta_n:n\ge  n_0\}$ of positive real numbers with $B^{\delta_n}\subseteq U$ for some $U\in \uc_n$ for all $n\ge  n_0$.  $X$ is said to be $\bs$-Lindel\"of \cite{cmk} if each $\bs$-cover contains a countable $\bs$-subcover.

Let $\{(X_n,d_n):n\in \nb\}$ be a family of metric spaces and let $\bk_n$ be a bornology on $X_n$ for each $n\in \nb$. Then the product bornology \cite{hh} on $\Pi_{n\in \nb}X_n$ has a base consisting of sets of the form $B=\Pi_{n\in \nb}B_n$ where $B_n\in \bk_n$ for all $n\in \nb$.

For two metric spaces $X$ and $Y$, $Y^X$ ($C(X, Y )$) stands for the set of all functions (continuous functions) from $X$ to $Y$. The commonly used
topologies on $C(X, Y )$ are the compact-open topology $\tau_k$, and the topology of pointwise convergence $\tau_p$. The corresponding spaces are, in general, respectively denoted by $(C(X, Y ), \tau_k)$ (resp. $C_k(X)$ when $Y = \rb$), and $(C(X, Y ), \tau_p)$ (resp. $C_p(X)$ when $Y = \rb$).\\

Let $\bk$ be a bornology with a closed base on $X$. Then the topology of strong uniform convergence $\tau_{\bk}^s$ is determined by a uniformity on $Y^X$ with a base consisting of all sets of the form
 \[ [B,\varepsilon]^s=\{(f,g): \exists \delta> 0 \,\, \text{for every} \,\, x\in B^{\delta},d(f(x),g(x))<\varepsilon \},\]
 for $B\in \bk, \varepsilon> 0$.

The topology of strong uniform convergence $\tau_{\bk}^s$ coincides with the topology of pointwise convergence $\tau_p$  if  $\bk= \fc$.

Throughout we use the convention that if $\bk$ is a bornology on $X$, then $X\not\in \bk$.

We end this section with some more basic observations about open $\bs$-covers.
First note that if $f:X\rightarrow Y$ is any map and if $\bk$ is a bornology on $X$, then the collection $f(\bk)=\{f(B):B\in \bk\}$ is a bornology on $f(X)$. Moreover, if $f$ is surjective then $f(\bk)$ is a bornology on $Y$. Note that if $\bk$ is a bornology on $X$ with a compact base $\bk_0$ and $f:X\rightarrow Y$ is a continuous function then $f(\bk)$ is a bornology on $f(X)$ with compact base $f(\bk_0)$.

\begin{Lemma}
\label{Lfc1}
Let $\bk$ be a bornology on $X$ with a compact base $\bk_0$ and $(Y,\rho)$ be a another metric space. Let $f:X\rightarrow Y$ be a continuous function on $X$. If $\uc$ be an open $f(\bk)^s$-cover ($\gamma_{f(\bk)^s}$-cover) of $f(X)$ then $\{f^{-1}(U):U\in \uc\}$ is an open $\bs$-cover ($\gamma_{\bk^s}$-cover) of $X$.
\end{Lemma}
\bp
Let $B\in \bk_0$. Since $\uc$ is an open $f(\bk)^s$-cover of $X$ so for $f(B)\in f(\bk)$ there is a $\varepsilon>0$ such that $f(B)^\varepsilon\subseteq U$ for some $U\in \uc$. As $f$ is continuous function on $B$ and $B$ is compact, $f$ is strongly uniformly continuous on $B$ \cite{bl} i.e. for $\varepsilon>0$ there is a $\delta>0$ such that $f(B^\delta)\subseteq f(B)^\varepsilon$. Therefore $f(B^\delta)\subseteq U$ i.e.  $B^\delta\subseteq f^{-1}(U)$. So $\{f^{-1}(U):U\in \uc\}$ is an open $\bs$-cover of $X$.
\ep

%\begin{Lemma}
%\label{Lfc2}
%Let $\bk$ be a bornology on $X$ with a compact base $\bk_0$ and $(Y,\rho)$ be a another metric space. Let $f:X\rightarrow Y$ be a continuous function on $X$.
%If $\uc$ is a ${f(\bk)^s}$-groupable cover of $f(X)$, then $\{f^{-1}(U):U\in \uc\}$ is a $\bs$-groupable cover of $X$.
%\end{Lemma}

\begin{Prop}
\label{Pf1}
Let $\bk$ be a bornology on $X$ with a compact base $\bk_0$ and $(Y,\rho)$ be a another metric space. Let $f:X\rightarrow Y$ be a continuous function on $X$.
Let $\Pi\in\{\S1,\Sf,\Uf\}$ , $\mathcal P, \mathcal Q  \in\{\oc,\Ga\}$.\\
If $X$ satisfies $\Pi(\mathcal P_\bs,\mathcal Q_\bs )$ then $f(X)$ satisfies $\Pi(\mathcal P_\fbs,\mathcal Q_\fbs)$.
\end{Prop}
\bp
We only show that if $X$ satisfies $\S1(\obs,\Gbs)$ then $f(X)$ satisfies $\S1(\oc_{f(\bk)^s},\Gamma_{f(\bk)^s})$. Before proceeding with the proof note that if $B$ is a compact and $U$ is an open subset of $X$ with $B\subseteq U$ then there is a $\delta>0$ such that $B^\delta\subseteq U$.

Let $\{\uc_n:n\in \nb\}$ be a sequence of open $f(\bk)^s$-covers of $f(X)$. By Lemma \ref{Lfc1}, $\uc'_n=\{f^{-1}(U):U\in \uc\}$ is an open $\bs$-cover of $X$ for each $n$. Apply $\S1(\obs,\Gbs)$ to the sequence $\{\uc'_n:n\in \nb\}$ to choose a $f^{-1}(U_n)\in \uc'_n$ for each $n$ such that $\{f^{-1}(U_n):n\in \nb\}$ is a $\gbs$-cover of $X$. We now show that the sequence $\{U_n:n\in \nb\}$ is a $\gamma_{f(\bk)^s}$-cover of $f(X)$. Let $B'\in f(\bk_0)$ and say $B' = f(B)$ where $B\in \bk_0$. Choose a $n_0\in \nb$ and a sequence $\{\delta_n:n\geq n_0\}$ of positive real numbers such that $B^{\delta_n}\subseteq f^{-1}(U_n)$ for $n\geq  n_0$ i.e. $f(B)\subseteq U_n$ for $n\geq  n_0$. Since $f(B)$ is compact, there is a $\varepsilon_n>0$ such that $f(B)^{\varepsilon_n}\subseteq U_n$ for each $n\geq n_0$. This shows that $\{U_n:n\in \nb\}$ is a $\gamma_{f(\bk)^s}$-cover of $f(X)$ and hence $f(X)$ satisfies $\S1(\oc_{f(\bk)^s},\Gamma_{f(\bk)^s})$.
\ep

\begin{Prop}
\label{Pf2}
Let $\bk$ be a bornology on $X$ with a compact base $\bk_0$. If $X$ satisfies $\S1(\Gbs,\Gbs)$ then every continuous image of $X$ into $\nb^\nb$ is bounded.
\end{Prop}
\bp
Let $\rho$ be the Baire metric on $\nb^\nb$ and $\varphi:X\rightarrow \nb^\nb$ be continuous. By Proposition \ref{Pf1}, $\varphi(X)$ satisfies $\S1(\Gamma_{\varphi(\bk)^s},\Gamma_{\varphi(\bk)^s})$. For $n,k\in \nb$, let $U^n_k=\{f\in \nb^\nb:f(n)\leq k\}$. Consider $\uc_n=\{U^n_k:k\in \nb\}$ for each $n$. Let $B\in \varphi(\bk_0)$. Using the compactness of $B$, we can easily find a $k_0$ and a sequence $\{\delta_k:k\geq k_0\}$ of positive real numbers such that such that $B^{\delta_k}\subseteq U^n_k$ for all $k\geq k_0$. Therefore $\uc_n$ is a $\gamma_{\varphi(\bk)^s}$-cover of $\varphi(X)$. Now apply $\S1(\Gamma_{\varphi(\bk)^s},\Gamma_{\varphi(\bk)^s})$ to $\{\uc_n:n\in \nb\}$ to choose a $U^n_{k_n}\in \uc_n$ for each $n$ such that $\{U^n_{k_n}:n\in \nb\}$ is an open $\gamma_{\varphi(\bk)^s}$-cover of $\varphi(X)$. Define a function $h:\nb\rightarrow \nb$ by $h(n)=k_n$ for $n\in \nb$. Let $f\in \varphi(X)$. Clearly there is a $n_0\in \nb$ such that $f\in U^n_{k_n}$ for all $n\geq n_0$ i.e. $f(n)\leq h(n)$ for all $n\geq n_0$ i.e. $f\leq^* h$. Hence $\varphi(X)$ is bounded in $\nb^\nb$.
\ep

\begin{Prop}
\label{Pf3}
Let $\bk$ be a bornology on $X$ with a compact base $\bk_0$. If $X$ satisfies $\Sf(\obs,\obs)$ then for every continuous function $\varphi:X\rightarrow \nb^\nb$, $\varphi(X)$ is not dominating.
\end{Prop}
\bp
By Proposition \ref{Pf1}, $\varphi(X)$  satisfies $\Sf(\oc_{\varphi(\bk)^s},\oc_{\varphi(\bk)^s})$. Consider $\uc_n=\{U^n_k:k\in \nb\}$ where $U^n_k=\{f\in \nb^\nb:f(n)\leq k\}$ for $n,k\in \nb$ which is a $\gamma_{\varphi(\bk)^s}$-cover of $\varphi(X)$. Apply
$\Sf(\oc_{\varphi(\bk)^s},\oc_{\varphi(\bk)^s})$ to $\{\uc_n:n\in \nb\}$ to choose a finite set $\vc_n\subseteq \uc_n$ for each $n$ such that $\cup_{n\in \nb}\vc_n$ is an open $\varphi(\bk)^s$-cover of $\varphi(X)$. Define a function $h:\nb\rightarrow \nb$ by $h(n)=\max\{k\in \nb:U^n_k\in \vc_n\}$ for $n\in \nb$. Now we show that for any $f\in \varphi(X)$, $f(n)\leq h(n)$ for infinitely many $n\in \nb$. Let $f\in \varphi(X)$. Choose a $B_0\in \varphi(\bk)$ such that $f\in B_0$. Since $\cup_{n\in \nb}\vc_n\in \oc_{\varphi(\bk)^s}$, in view of \cite[Proposition 3.1]{dcpdsd}, there are $\delta_n>0$ such that $B_0^{\delta_n}\subseteq U^n_k$ for some $U^n_k\in \vc_n$ for infinitely many $n\in \nb$ i.e. $f\in U^n_k$ for some $U^n_k\in \vc_n$ for infinitely many $n\in \nb$ i.e. $f(n)\leq h(n)$ for infinitely many $n\in \nb$. Hence $\varphi(X)$ is not dominating.
\ep

%\textbf{Problem}: Is the converse of Proposition \ref{Pf2} and \ref{Pf3} true?

\section{Results related to cardinality and Ramsey Theory}

\subsection{Results concerning cardinality}

Recall that $X$ is called a $\gbs$-set \cite{dcpdsd} if $X$ satisfies $\S1(\obs,\Gbs)$. An equivalent condition is that every open $\bs$-cover of $X$ contains a countable set which is a $\gbs$-cover of $X$ (see \cite{cmk}).
\begin{Th}
\label{Tcard3}
Let $\bk$ be a bornology on $X$ with a closed base $\bk_0$ and $X$ be $\bs$-Lindel\"of. If $|\bk_0|<\pf$ then $X$ is a $\gbs$-set.
\end{Th}
\bp
Let $\uc$ be an open $\bs$-cover of $X$. Enumerate $\uc$ bijectively as $\{U_n:n\in \nb\}$. By \cite[Proposition 3.1]{dcpdsd}, for $B\in \bk$ there are $\delta_n>0$ and $U_n\in \uc$ such that $B^{\delta_n}\subseteq U_n$ for infinitely many $n$. Define $A_B=\{n\in \nb:B^{\delta_n}\subseteq U_n\}$ for $B\in \bk$. Clearly each $A_B$ is an infinite subset of $\mathbb{N}$ and now consider the family $\ac=\{A_B:B\in \bk_0\}$.  Let $\{A_{B_1},\dots,A_{B_k}\}$ be any finite subfamily of $\ac$. Since $B_1\cup\dots \cup B_k\in \bk$, there are  $\delta_n>0$ and $U_n\in \uc$ such that $(B_1\cup\dots \cup B_k)^{\delta_n}\subseteq U_n$ for infinitely many $n$. This means that $A_{B_1}\cap\dotsc \cap A_{B_k}$ must be infinite and consequently we can conclude that any finite subfamily of $\ac$ has infinite intersection. Now in view of our assumption that $|\bk_0|<\pf$, we can choose an infinite subset $P$ of $\nb$ such that for each $B\in \bk_0$, $P\setminus A_B$ is finite. Enumerating $P$ as $\{n_k:k\in \nb\}$, set $\vc=\{U_{n_k}\in \uc:k\in \nb\}$. We show that $\vc$ is a $\gbs$-cover. Let $B\in \bk_0$. Using the fact that $P\setminus A_B$ is finite and the definition of $A_B$, there is a $k_0\in \nb$ and a sequence $\{\delta_{n_k}:k\geq k_0\}$ of positive real numbers such that $B^{\delta_{n_k}}\subseteq U_{n_k}$ for all $k\geq k_0$. This shows that $\vc\subseteq \uc$ is a $\gbs$-cover of $X$ and hence $X$ is a $\gbs$-set.
\ep

\begin{Th}
\label{Tcard4}
Let $\bk$ be a bornology on $X$ with a closed base $\bk_0$ and $X$ be $\bs$-Lindel\"of. If $|\bk_0|<\cov(\mc)$ then $X$ satisfies $\S1(\obs,\obs)$.
\end{Th}
\bp
Let $\{\uc_n:n\in \nb\}$ be a sequence of open $\bs$-covers of $X$. Enumerate each $\uc_n$ bijectively as $\{U^n_m:m\in \nb\}$. For $B\in \bk_0$, choose a $\delta>0$ and a $U^n_m\in \uc_n$ such that $B^\delta\subseteq U^n_m$. Define a function $f_B:\nb\rightarrow \nb$ by $f_B(n)=\min\{m\in \nb:B^\delta\subseteq U^n_m$ for some $\delta>0\}$ for each $n\in \nb$. Consider the set $\{f_B:B\in \bk_0\}$. Since $|\bk_0|<\cov(\mc)$, there is a $f:\nb\rightarrow \nb$ such that for each $B\in \bk_0$, $\{n\in \nb:f_B(n)=f(n)\}$ is infinite \cite[Theorem 5]{fm} (for our consideration the fact that this particular set is non-empty is the most important fact). We show that $\{U^n_{f(n)}:n\in \nb\}$ is the required open $\bs$-cover. Let $B\in \bk$. There is a $B_0\in \bk_0$ with $B \subseteq B_0$. By definition of $f_{B_0}$, $B_0^\delta \subseteq U^n_{f_{B_0}(n)}$ for some $\delta>0$. Since the set $\{n:f_{B_0}(n)=f(n)\} \neq \emptyset$, choosing an appropriate $k\in \nb$ so that $f_{B_0}(k)=f(k)$ we observe that $B^\delta \subseteq U^k_{f(k)}$ for some $\delta>0$. So, as claimed, $\{U^n_{f(n)}:n\in \nb\}$ is an open $\bs$-cover of $X$ showing that $X$ satisfies $\S1(\obs,\obs)$.
\ep

\begin{Th}
\label{Tcard2}
 Let $\bk$ be a bornology on $X$ with a closed base $\bk_0$ and $X$ be $\bs$-Lindel\"of. If $|\bk_0|<\df$ then $X$ satisfies $\Sf(\obs,\obs)$.
\end{Th}

\bp
Let $\{\uc_n:n\in \nb\}$ be a sequence of open $\bs$-covers of $X$. Enumerate each $\uc_n$ bijectively as $\uc_n=\{U^n_m:m\in \nb\}$ for each $n$. For $B\in\bk_0$, there are $\delta>0$ and  $U^n_m\in \uc_n$ such that $B^\delta\subseteq U^n_m$. Define a function $f_B:\nb\rightarrow \nb$ by $f_B(n)=\min\{m\in \nb:B^\delta\subseteq U^n_m$ for some $\delta>0\}$. Consider the set $\{f_B:B\in \bk_0\}$. Since $|\bk_0|<\df$,  the family $\{f_B:B\in \bk_0\}$ is not dominating. Therefore there is a function $g:\nb\rightarrow \nb$ such that for any $B\in \bk_0$, $f_B(n)<g(n)$ for infinitely many $n$. Define $\vc_n=\{U^n_m:m\le  g(n)\}$ for each $n$. Clearly $\vc_n$ is a finite subset of $\uc_n$ for each $n$ and $\cup_{n\in \nb}\vc_n$ is an open $\bs$-cover of $X$. Hence $X$ satisfies $\Sf(\obs,\obs)$.
\ep

\begin{Th}
\label{Tcard1}
Let $\bk$ be a bornology on $X$ with a closed base $\bk_0$. If $|\bk_0|<\bb$ then $X$ satisfies $\S1(\Gbs,\Gbs)$.
\end{Th}
\bp
Let $\{\uc_n:n\in \nb\}$ be a sequence of $\gbs$-covers of $X$, where $\uc_n=\{U^n_m:m\in \nb\}$ for each $n$. For each $B\in\bk_0$, there exist a $m_0\in \nb$ and a sequence $\{\delta_m:m\ge m_0\}$ of positive real numbers such that $B^{\delta_m}\subseteq U^n_m$ for $m\geq  m_0$. Define a function $f_B:\nb\rightarrow \nb$ by $f_B(n)=\min\{k\in \nb:$ for all $m\ge k,B^{\delta_m}\subseteq U^n_m\}$. Consider the set $\{f_B:B\in \bk_0\}$. Since $|\bk_0|<\bb$, there is a $g:\nb\rightarrow \nb$ such that $f_B\leq ^*g$ for all $B\in \bk_0$. We now show that $\{U^n_{g(n)}:n\in \nb\}$ is a $\gbs$-cover of $X$. Let $B\in \bk$ and choose a $B_0\in \bk_0$ such that $B\subseteq B_0$. As $f_{B_0}\leq ^* g$, there is a $n_1\in \nb$ such that $f_{B_0}(n)\leq  g(n)$ for all $n\geq  n_1$. From the definition of $f_{B_0}$, one can construct a sequence $\{\varepsilon_n:n\ge n_1\}$ of positive real numbers such that $B_0^{\varepsilon_n}\subseteq U^n_{g(n)}$ for all $n\geq  n_1$ and so $B^{\varepsilon_n}\subseteq B_0^{\varepsilon_n}\subseteq U^n_{g(n)}$ for all $n\geq  n_1$. Therefore $\{U^n_{g(n)}:n\in \nb\}$ is a $\gbs$-cover of $X$ which shows that $X$ satisfies $\S1(\Gbs,\Gbs)$.
\ep

\subsection{Ramsey theoretic results}
%This section concerns with partition relations.
We now formulate $\S1$ and  $\Sf$-type selection principles using certain partition relations involving bornological covers.
%We show that the selection principles of $\S1$,$\Sf$-type involving the family $\obs$, $\Gbs$ of bornological covers can be characterized in terms of certain partition relations.
The following observation about open $\bs$-covers will be useful in this context.
\begin{Lemma}
\label{Lp1}
Let $\bk$ be a bornology on $X$ with closed base. If an open $\bs$-cover $\uc$ can be expressed as a union of finite number of subfamilies of $\uc$, then at least one of them must be an open $\bs$-cover of $X$.
\end{Lemma}
\bp
Let $\uc$ be an open $\bs$-cover and $\uc=\displaystyle{\cup_{i=1}^k} \uc_i$. If on the contrary, no $\uc_i$ is an open $\bs$-cover of $X$, then for each $\uc_i$ there is a $B_i\in \bk$ such that for any $\delta>0$ and for any $U\in \uc_i$, $B^\delta_i\nsubseteq U$. Take $B=B_1\cup \dots \cup B_k\in \bk$. It is clear that for any $\delta>0$, $B^\delta\nsubseteq U$ for any $U\in \uc$. This contradicts that $\uc$ is an open $\bs$-cover of $X$. Hence at least one of $\uc_i$'s must be an open $\bs$-cover of $X$.
\ep

The method of proof of the following theorem is adapted from \cite[Theorem 10]{cooc1} and \cite[Theorem 6.2]{cooc2} with suitable modifications.
\begin{Th}
\label{Tp1}
Let $\bk$ be a bornology on $X$ with closed base. The following statements are equivalent.

$(1)$ $X$ satisfies $\Sf(\obs,\obs)$.

$(2)$ $X$ satisfies $\obs\rightarrow \lceil \obs\rceil^2_2$, provided $X$ is $\bs$-Lindel\"{o}f.
\end{Th}
\bp
$(1)\Rightarrow (2)$.
Let $\uc=\{U_n:n\in \nb\}$ be an open $\bs$-cover of $X$ and $f:[\uc]^2\rightarrow \{1,2\}$ be a coloring. Clearly $\uc$ can be expressed as $\tc_1\cup \tc_2$, where $\tc_1=\{V\in \uc:f(\{U_1,V\})=1\}$ and $\tc_2=\{V\in \uc:f(\{U_1,V\})=2\}$. By Lemma \ref{Lp1}, at least one of $\tc_1$ and $\tc_2$ must be an open $\bs$-cover of $X$. Let $i_1\in \{1,2\}$ be such that $\uc_1=\tc_{i_1}$ is an open $\bs$-cover. Inductively, choose sequences $\{\uc_n:n\in \nb\}$ and $\{i_n:n\in \nb\}$ such that $\uc_n=\{V\in \uc_{n-1}:f({U_n,V})=i_n\}$ is an open $\bs$-cover for $n>1$. Now apply $\Sf(\obs,\obs)$ to $\{\uc_n:n\in \nb\}$ to obtain a sequence $\{\vc_n:n\in \nb\}$ where $\vc_n$ is a finite subset of $\uc_n$ for each $n$ such that ${\cup_{n\in \nb}}\vc_n\in \obs$. Further we can assume that $\vc_n$'s are pairwise disjoint and there is a $i\in \{1,2\}$ such that for each $U_m\in {\cup_{n\in \nb}}\vc_n$, $i_m=i$. Choose $k_1\in \nb$ large enough so that $i\leq k_1$ whenever $U_i\in \vc_1$. Further choose $k_2>k_1$ so that $i\leq k_2$ whenever $U_i\in \vc_2$ and so on. Again choose a sequence $l_1<l_2< \cdots$ such that ${\cup_{j>l_1}} \vc_j\subseteq \uc_{k_1}$ and ${\cup_{j>l_m}} \vc_j\subseteq \uc_{k_{l_m}}$ for $m>1$. Take $\zc_n={\cup_{l_n\leq j<l_{n+1}}}\vc_j$ for each $n\in \nb$. Clearly ${\cup_{n\in \nb}}\zc_n$ is an open $\bs$-cover and $\{\zc_n:n\in \nb\}$ is a partition of it into pairwise disjoint finite sets. Again by Lemma \ref{Lp1}, at least one of ${\cup_{n\in \nb}}\zc_{2n}$ and  ${\cup_{n\in \nb}}\zc_{2n-1}$ must be an open $\bs$-cover and without any loss of generality assume that ${\cup_{n\in \nb}}\zc_{2n}$ is the one. Write $\vc={\cup_{n\in \nb}}\zc_{2n}$. Define $\varphi:\vc\rightarrow \nb$ by $\varphi(V)=n$ if $V\in \zc_{2n}$. Then for all $V,W\in \vc$ with $\varphi(V)\neq \varphi(W)$, $f(\{V,W\})=i$.

$(2)\Rightarrow(1)$. Let $\{\uc_n:n\in \nb\}$ be a sequence of open $\bs$-covers of $X$ and let $\uc_n=\{U^n_m:m\in \nb\}$ for each $n\in \nb$. Consider the collection $\vc=\{U^1_m\cap U^m_k:m,k\in \nb\}$. To see that $\vc$ is an open $\bs$-cover, let $B\in \bk$. There are $\delta_1, \delta_2>0$ and $U^1_m\in \uc_1$, $U^m_k\in \uc_m$ such that $B^{\delta_1}\subseteq U^1_m$ and $B^{\delta_2}\subseteq U^m_k$. Choosing $\delta=\min\{\delta_1, \delta_2\}$, we get that $B^\delta\subseteq U^1_m\cap U^m_k$. Define $f:[\vc]^2\rightarrow \{1,2\}$ by
\[f(\{U^1_{m_1}\cap U^{m_1}_{k_1}, U^1_{m_2}\cap U^{m_2}_{k_2}\})=\\
\begin{cases}
1&\:\: \text{if}\:\: m_1=m_2\\

2&\:\: \text{if}\:\: m_1\neq m_2
\end{cases}\].

By $(2)$, there are a $\bs$-subcover $\wc=\{V_j:j\in \nb\}$ of $\vc$ where $V_j=U^1_{m_j}\cap U^{m_j}_{k_j}$, a finite-to-one function $\varphi:\wc\rightarrow \nb$ and a $i\in \{1,2\}$ such that for all $V_j,V_p\in \wc$ whenever $\varphi(V_j)\neq \varphi(V_p)$, $f(\{V_j,V_p\})=i$.

If $i=1$. Then for all $j,p\in \nb$, $m_j=m_p$ i.e. $m_j=m_1$ for all $j\in \nb$ which  implies that every element of $\wc$ refines $U^1_{m_1}$. But this contradicts that $\wc$ is an open $\bs$-cover of $X$. Therefore we must have $i=2$.

Let $\zc=\{U^{m_j}_{k_j}:j\in \nb\}$. Clearly $\zc$ is an open $\bs$-cover as $\wc$ refines $\zc$. Choose $\zc_n=\{U^{m_j}_{k_j}:m_j=n\}$. $\zc_n$'s are either empty or finite. Otherwise there will be infinitely many $j\in \nb$ with $m_j=n$. Since $\varphi$ is finite-to-one, choose $j,p$ with $m_j=m_p=n$ and $\varphi(V_j)\neq \varphi(V_p)$, then $f(\{V_j,V_p\})=2$ i.e. $m_j\neq m_p$, which is a contradiction. Hence $\{\zc_n:n\in \nb\}$ is a sequence of finite sets with $\zc_n\subseteq \uc_n$ for each $n\in \nb$ witnessing $\Sf(\obs,\obs)$.
\ep

Using the method of proof of  Theorem \ref{Tp1}, the induction argument on $n$ and $k$ and the usual method for proving the Ramsey theoretic statements for $n > 2$, $k > 2$ (see, for example, \cite[Theorem 1]{coocvii} and \cite[Theorem 25]{cooc1}), we obtain the following result.

\begin{Th}
\label{Tp2}
Let $\bk$ be a bornology on $X$ with closed base. The  following statements are equivalent.

$(1)$ $X$ satisfies $\S1(\obs,\obs)$.

$(2)$ $X$ satisfies $\obs\rightarrow \ (\obs)^2_2$, provided $X$ is $\bs$-Lindel\"{o}f.

$(3)$ $X$ satisfies $\obs\rightarrow \ (\obs)^n_k$ for each $n,k\in \nb$, provided $X$ is $\bs$-Lindel\"{o}f.
\end{Th}

\begin{Th}
\label{Tp3}
Let $\bk$ be a bornology on $X$ with closed base. The  following statements are equivalent.

$(1)$ $X$ satisfies $\S1(\obs,\Gbs)$.

$(2)$ ONE does not have a winning strategy in $\G1(\obs,\Gbs)$ on $X$.

$(3)$ $X$ satisfies $\obs\rightarrow \ (\Gbs)^n_k$ for each $n,k\in \nb$, provided $X$ is $\bs$-Lindel\"{o}f.

%$(4)$ Every open $\bs$-cover of $X$ contains a countable set which is a $\gbs$-cover of $X$.
\end{Th}
\bp
The equivalence $(1)\Leftrightarrow (2)$ is already proved in \cite[Theorem 3.7]{dcpdsd}.
%and the implication $(4)\Rightarrow (1)$ is due to \cite[Theorem 2.8]{cmk}.
We prove the implications $(2)\Rightarrow (3)$ and $(3)\Rightarrow (1)$.

$(2)\Rightarrow (3)$. For this proof we follow the line of argument of \cite[Theorem 30]{cooc1}. Let $\uc$ be an open $\bs$-cover of $X$ and let $f:[\uc]^2\rightarrow \{1,\dots,k\}$ be a coloring. Enumerate $\uc$ bijectively as $\{U_{(n)}:n\in \nb\}$. We define a strategy $\sigma$ for ONE in $\G1(\obs,\Gbs)$ as follows.

Define $\sigma(\emptyset)=\uc$, the first move of ONE. TWO responds by selecting $U_{(n_1)}\in \uc$. Now $\uc\setminus \{U_{(n_1)}\}$ can be expressed as ${\cup_{i=1}^k} \tc_i$, where $\tc_i=\{V\in \uc:f(\{U_{(n_1)},V\})=i\}$ for $i=1,\dotsc,k$. By Lemma \ref{Lp1}, at least one of $\tc_i$'s is an open $\bs$-cover of $X$ and assume it to be $\tc_{i_{(n_1)}}$. Clearly $\tc_{i_{(n_1)}}=\{V\in \sigma(\emptyset):f(\{U_{(n_1)},V\})=i_{(n_1)}\}$. Enumerate $\tc_{i_{(n_1)}}$ bijectively as $\{U_{(n_1,m)}:m\in \nb\}$. Now define $\sigma(U_{(n_1)})=\{U_{(n_1,m)}:m\in \nb\}$. TWO responds by selecting $U_{(n_1,n_2)}$. Continuing this way, we obtain a $i_{(n_1,\dotsc,n_r)}\in \{1,\dotsc,k\}$ and an open $\bs$-cover $\{V\in \tau(U_{(n_1)},\dotsc,U_{(n_1,\dotsc,n_{r-1})})
:f(\{U_{(n_1,\dotsc,n_r)},V\})=i_{(n_1,\dotsc,n_r)}\}$ which is enumerated bijectively as $\{U_{(n_1,\dotsc,n_r,m)}:m\in \nb\}$. Then define $\sigma(U_{(n_1)},\dotsc,U_{(n_1,\dotsc,n_r)})
=\{U_{(n_1,\dotsc,n_r,m)}:m\in \nb\}$. TWO responds by selecting $U_{(n_1,\dotsc,n_{r+1})}$ and so on.

This defines  a strategy $\sigma$ for ONE in $\G1(\obs,\Gbs)$. By $(2)$, $\sigma$ is not a winning strategy for ONE. Consider a play
\[\sigma(\emptyset), U_{(n_1)},\dotsc, \sigma(U_{(n_1)},\dotsc,U_{(n_1,\dotsc,n_r)}), U_{(n_1,\dotsc,n_{r+1})},\dotsc~\]
which is lost by ONE. Therefore
\[U_{(n_1)},\dotsc, U_{(n_1,\dotsc,n_{r+1})},\dotsc\] should form a $\gbs$-cover of $X$. Since for all $r\in \nb$, $i_{(n_1,\dotsc,n_r)}\in \{1,\dotsc,k\}$, there is a $i\in \{1,\dotsc,k\}$ such that for infinitely many $r\in \nb$, $i_{(n_1,\dotsc,n_r)}=i$. Let $\vc=\{U_{(n_1,\dotsc,n_r)}:i_{(n_1,\dotsc,n_r)}=i\}$. As every infinite subset of a $\gbs$-cover is a $\gbs$-cover, $\vc$ is a $\gbs$-cover of $X$. Finally note that, from the construction of the game, for any $r,s\in \nb$ with $r<s$ we must have  $f(\{U_{(n_1,\dotsc,n_r)},U_{(n_1,\dotsc,n_s)}\})
=i_{(n_1,\dotsc,n_r)}$. Therefore for any $V,W\in \vc$ with $V\neq W$, $f(\{V,W\})=i$. Hence $(3)$ holds.

%$(3)\Rightarrow (4)$. Let $\uc$ be an open $\bs$-cover of $X$. Let $\{\uc_n:n\in \nb\}$ be a partition of $\uc$ into nonempty finite sets. Define a coloring $f:[\uc]^2\rightarrow \{1,2\}$ by $f(\{U,V\})=1$ if $U,V\in \uc_n$ for some $n\in \nb$ and $f(\{U,V\})=2$ otherwise.
%%We show that $f$ is a coloring. Let $U_1\in \uc$ and $\vc$ be a $\bs$-subcover of $\uc$. Now we can write $\uc=\{V\in \uc:f(\{U_1,V\})=1\}\cup \{V\in \uc:f(\{U_1,V\})=2\}$. By Lemma \ref{Lp1},
%By $(3)$, there is a $\vc\subseteq \uc$ which is a $\gbs$-cover of $X$. Hence $(4)$ holds.

$(3)\Rightarrow (1)$. Let $\{\uc_n:n\in \nb\}$ be a sequence of open $\bs$-covers of $X$. Let $\uc_n=\{U^n_p:p\in \nb\}$. Consider the collection $\vc_n=\{U^1_{p_1}\cap U^2_{p_2}\cap \dots \cap U^n_{p_n}: U^i_{p_i}\in \uc_i, 1<p_1<p_2<\cdots<p_n\}$ for each $n\in \nb$. Then $\vc_n$'s are open $\bs$-covers of $X$. For convenience, we write $\vc_n=\{V^n_m:m\in \nb\}$ for $n\in \nb$. Now define $\vc=\{V^1_m\cap V^m_k:m,k\in \nb\}$. Clearly $\vc$ is an open $\bs$-cover of $X$. Define a coloring $f:[\vc]^2\rightarrow \{1,2\}$ by
\[f(\{V^1_{m_1}\cap V^{m_1}_{k_1},V^1_{m_2}\cap V^{m_2}_{k_2}\})=\\
\begin{cases}
1& \:\: \text{if}\:\: m_1=m_2\\

2& \:\: \text{if}\:\: m_1\neq m_2
\end{cases}\]
By $(3)$, there are a $\wc\in \Gbs$ with $\wc\subseteq \vc$, a function $\varphi:\wc\rightarrow \nb$ and a $i\in \{1,2\}$ such that for $U,V\in \wc$, $f(\{U,V\})=i$.

Clearly $i=2$. As $i=1$ would imply that every element of $\wc$ refines $V^1_{m_1}$, contradicting that $\wc$ is a $\gbs$-cover of $X$.
%%%%%%%%%%%%%%%%%%%%%%%%%%%%%%%%%%%%%%%%%%%%%%%%
%%%%%change required in the following%%%%%%%%%%%
%%%%%%%%%%%%%%%%%%%%%%%%%%%%%%%%%%%%%%%%%%%%%%%%

Enumerate $\wc$ as $\{V^1_{m_j}\cap V^{m_j}_{k_j}:j\in \nb\}$ where $m_1<m_2<\cdots$. Consider $\{ V^{m_j}_{k_j}:j\in \nb\}$ which is a $\gbs$-cover of $X$ as $\wc$ is a $\gbs$-cover of $X$.  The elements of $\vc_n$ are of the form $U^1_{p_1}\cap U^2_{p_2}\cap \dots \cap U^n_{p_n}$ where $U^i_{p_i}\in \uc_i$.
For each $n\in [1,m_1]$, we choose the $n$-th coordinate $U^n_{p_n}$ in the representation of $V^{m_1}_{k_1}$ and for each $j>1$ and $n\in (m_{j-1},m_j]$, choose the $n$th coordinate $U^n_{p_n}$ in the representation of $V^{m_j}_{k_j}$. Clearly $\{U^n_{p_n}:n\in \nb\}$ is a $\gbs$-cover of $X$ and hence $X$ satisfies $\S1(\obs,\Gbs)$.
\ep

%%%%%%%%%%%%%%%%%%%%%%%%%%%%%%%%%%%%%%%%%%%%%
%%%%%%%%%%%%%%%%%%%%%%%%%%%%%%%%%%%%%%%%%%%%%

\section{The $\bs$-Hurewicz property and the $\bs$-Gerlits-Nagy property}

\subsection{Some observations on $\bs$-Hurewicz property}

%\vspace{0.5 cm}

In \cite{dcpdsd}, the notion of strong $\bk$-Hurewicz property (or in short, $\bs$-Hurewicz property) was introduced and several of its basic properties were established. In this section first we again look back at this very important property and present some more new observations.

 Recall that $X$ is said to have the $\bs$-Hurewicz property  if for each sequence $\{\uc_n:n\in \nb\}$ of open ${\bk^s}$-covers of $X$, there is a sequence  $\{\vc_n:n\in \nb\}$ where $\vc_n$ is a finite subset of $\uc_n$ for each $n\in \nb$, such that for every $B\in \bk$ there exist a $n_0\in \nb$ and a sequence $\{\delta_n:n\ge  n_0\}$ of positive real numbers satisfying $B^{\delta_n}\subseteq U$ for some $U\in \vc_n$ for all $n\ge  n_0$.

\begin{Prop}
\label{Pfc1}
Let $\bk$ be a bornology on $X$ with a compact base $\bk_0$ and $(Y,\rho)$ be a another metric space. Let $f:X\rightarrow Y$ be a continuous function on $X$. If $X$ has the $\bs$-Hurewicz property then $f(X)$ has the $f(\bk)^s$-Hurewicz property.
\end{Prop}
\bp
Let $\{\uc_n:n\in \nb\}$ be a sequence of open $f(\bk)^s$-covers of $f(X)$. By Lemma \ref{Lfc1}, for each $n\in \nb$, $\uc_n'=\{f^{-1}(U):U\in \uc_n\}$ is an open $\bs$-cover of $X$. In view of the $\bs$-Hurewicz property of $X$, we can find a sequence $\{\vc_n':n\in \nb\}$ of finite sets with $\vc_n'\subseteq \uc_n'$ for each $n$, such that for each $B\in \bk$ there exist a $n_0$ and a sequence $\{\delta_n:n\ge n_0\}$ of positive reals satisfying $B^{\delta_n}\subseteq f^{-1}(U)$ for some $f^{-1}(U)\in \vc_n'$ for all $n\ge n_0$. For each $n$ choose $\vc_n=\{U\in \uc_n:f^{-1}(U)\in \vc_n'\}$. We show that $\{\vc_n:n\in \nb\}$ witnesses the $f(\bk)^s$-Hurewicz property for $f(X)$. Let $B'\in f(\bk_0)$ and say $B' = f(B)$ where $B\in \bk_0$. We can find a $n_0\in \nb$ and a sequence $\{\delta_n:n\ge n_0\}$ of positive reals satisfying $B^{\delta_n}\subseteq f^{-1}(U)$ for some $f^{-1}(U)\in \vc_n'$ for all $n\ge n_0$. Subsequently $f(B)\subseteq U$ for some $U\in \vc_n$ for all $n\ge n_0$. Since $f(B)$ is compact, for each $n\geq n_0$ there is a $\varepsilon_n>0$ such that $f(B)^{\varepsilon_n}\subseteq U$. This shows that $\{\vc_n:n\in \nb\}$ witnesses the $f(\bk)^s$-Hurewicz property for $f(X)$.
\ep

\begin{Prop}
\label{Pfc2}
Let $\bk$ be a bornology on $X$ with a compact base $\bk_0$. If $X$ has the $\bs$-Hurewicz property then every continuous image of $X$ into $\nb^\nb$ is bounded.
\end{Prop}
\bp
Let $\rho$ be the Baire metric on $\nb^\nb$ and $\varphi:X\rightarrow \nb^\nb$ be continuous. By Proposition \ref{Pfc1}, $\varphi(X)$ has  the $\varphi(\bk)^s$-Hurewicz property. Consider $\uc_n=\{U^n_k:k\in \nb\}$ where $U^n_k=\{f\in \nb^\nb:f(n)\leq k\}$ for $n,k\in \nb$ which is an open ${\varphi(\bk)^s}$-cover of $\varphi(X)$.
%Let $B\in \varphi(\bk_0)$. Choosing $\delta=\frac{1}{2(n+1)}$ and using compactness of $B$, one can easily show that $B^\delta\subseteq U^n_k$ for some $U^n_k\in \uc_n$.
%$B=\cup_{f\in B} S_\rho(f,\frac{1}{2(n+1)})$, since $B$ is compact there are $f_1,f_2,\dots, f_r$ such that $B=S_\rho(f_1,\frac{1}{2(n+1)})\cup \dots \cup S_\rho(f_r,\frac{1}{2(n+1)})$. Choose $k=\max\{f_1(n),\dots,f_r(n)\}$. Then $\{f_1,\dots,f_r\}\subseteq U^n_k$. It is easily seen that $S_\rho(f_i,\frac{1}{n+1})\subseteq U^n_k$ for $i=1,\dots,r$. Choose $\delta=\frac{1}{2(n+1)}$. We show that $B^\delta\subseteq U^n_k$. Let $g\in B^\delta$. Clearly $\rho(g,f)<\delta$ for some $f\in B$. Now for $f\in B$, $\rho(f,f_i)<\frac{1}{2(n+1)}$ for some $i\in \{1,\dots,r\}$. So $\rho(g,f_i)<\frac{1}{n+1}$ i.e. $g\in S_\rho(f_i,\frac{1}{n+1})\subseteq U^n_k$. Therefore $B^\delta\subseteq U^n_k$.

Apply $\varphi(\bk)^s$-Hurewicz property to $\{\uc_n:n\in \nb\}$ to obtain $\{\vc_n:n\in \nb\}$ where each $\vc_n$ is a finite subset of $\uc_n$ such that for each $B\in \varphi(\bk_0)$ there exist a $n_0$ and a sequence $\{\delta_n:n\geq n_0\}$ of positive real numbers satisfying $B^{\delta_n}\subseteq U^n_k$ for some $U^n_k\in \vc_n$ and for all $n\geq n_0$. Define a function $h:\nb\rightarrow \nb$ by $h(n)=\max\{k:U^n_k\in \vc_n\}$. We now show that for any $f\in \varphi(X)$, $f\leq^* h$ holds. For $f\in \varphi(X)$ choose a $B_0\in\varphi(\bk_0)$ such that $f\in B_0$. Now choose a $n_0$ and a sequence $\{\delta_n:n\geq n_0\}$ of positive real numbers satisfying $B_0^{\delta_n}\subseteq U^n_k$ for some $U^n_k\in \vc_n$ and for all $n\geq n_0$ i.e. $f\in U^n_k$ for some $U^n_k\in \vc_n$ and for all $n\geq n_0$ i.e. $f(n)\leq h(n)$ for all $n\geq n_0$. Thus $f\leq^* h$ holds and hence $\varphi(X)$ is bounded.
\ep

%\textbf{Problem}: Is the converse of Proposition \ref{Pfc2} true?\\

The following result provides an important insight into the spaces with $\bs$-Hurewicz property in terms of the cardinality of the base of the bornology and the significance of the result is due to the fact that several known bornolocical spaces have closed bases.
\begin{Th}
\label{Thb}
Let $\bk$ be a bornology on $X$ with a closed base $\bk_0$. %The following statements are equivalent.
%$(1)$ $|\bk_0|<\bb$.
%$(2)$ $X$ has the $\bs$-Hurewicz property.
If $|\bk_0|<\bb$ then $X$ has the $\bs$-Hurewicz property.
\end{Th}
\bp
%$(1)\Rightarrow (2)$.
Let $\{\uc_n:n\in \nb\}$ be a sequence of open $\bs$-covers of $X$. Let $\uc_n=\{U^n_m:m\in \nb\}$ for $n\in \nb$. For $B\in \bk_0$ there are $\delta>0$ and $U^n_m\in \uc_n$ such that $B^\delta\subseteq U^n_m$. For each $B\in \bk_0$, define a function $f_B:\nb\rightarrow \nb$ by $f_B(n)=\min\{m\in \nb:B^\delta\subseteq U^n_m$ for some $\delta>0\}$.
Consider the set $\{f_B:B\in \bk_0\}$. Since $|\bk_0|<\bb$, there is a $g:\nb\rightarrow \nb$ such that $f_B\leq ^*g$ for all $B\in \bk_0$. For each $n\in \nb$, choose $\vc_n=\{U^n_m:m\leq  g(n)\}$. We claim that $\{\vc_n:n\in \nb\}$ witnesses the $\bs$-Hurewicz property.

To see this we will show that for $B\in \bk$ there is a $n_0\in \nb$ and a sequence $\{\delta_n:n\geq  n_0\}$ of positive real numbers such that $B^{\delta_n}\subseteq U$ for some $U\in \vc_n$ for all $n\geq n_0$. For $B\in \bk$ choose $B_0\in \bk_0$ with $B\subseteq B_0$. Since $f_{B_0}\leq ^*g$, there is a $n_0\in \nb$ such that $f_{B_0}(n)\leq  g(n)$ for all $n\geq  n_0$. Now by definition of $f_{B_0}$, for each $n\geq  n_0$ there is a $\delta_n>0$ such that $B_0^{\delta_n}\subseteq U^n_{f_{B_0}(n)}$. Clearly $U^n_{f_{B_0}(n)}\in \vc_n$ for all $n\geq  n_0$. So we have a sequence $\{\delta_n:n\geq  n_0\}$ of positive real numbers such that $B^{\delta_n}\subseteq B_0^{\delta_n}\subseteq U$ for some $U\in \vc_n$ for all $n\geq  n_0$. Hence $X$ has the $\bs$-Hurewicz property.

%$(2)\Rightarrow (1)$. On the contrary, let us assume that $(1)$ is not true i.e. $|\bk_0|\geq  \bb$. Let $\{\uc_n:n\in \nb\}$ be a sequence of open $\bs$-covers of $X$. As before write $\uc_n=\{U^n_m:m\in \nb\}$. For $B\in \bk_0$ there are $\delta>0$ and $U^n_m\in \uc_n$ such that $B^\delta\subseteq U^n_m$. For each $B\in \bk_0$, consider the function $f_B:\nb\rightarrow \nb$ defined above.
%As $|\bk_0|\geq  \bb$, for the family $\{f_B:B\in \bk_0\}$ and for any $g:\nb\rightarrow \nb$ there must be a $B\in \bk_0$ such that $g(n)<f_B(n)$ for infinitely many $n\in \nb$. By $(2)$, there is a sequence $\{\vc_n:n\in \nb\}$ of finite sets with $\vc_n\subseteq \uc_n$ for each $n\in \nb$ witnessing the $\bs$-Hurewicz property of $X$. Now define a function $g:\nb\rightarrow \nb$ by $g(n)=\max\{m\in \nb:U^n_m\in \vc_n\}$. For this particular $g$, one should find a $B_0\in \bk_0$ such that $g(n)<f_{B_0}(n)$ for infinitely many $n\in \nb$. Now from the definition of $f_{B_0}$, it follows that, for any $\delta>0$ $B_0^\delta\nsubseteq U^n_m$ whenever $m\leq  g(n)$ for infinitely many $n\in \nb$. But then for any sequence $\{\delta_n:n\in \nb\}$ of positive real numbers and for any $U^n_m\in \vc_n$, $B_0^{\delta_n}\nsubseteq U^n_m$  for infinitely many $n\in \nb$ which contradicts that $\{\vc_n:n\in \nb\}$ witnesses the $\bs$-Hurewicz property. Therefore we must have $|\bk_0|<\bb$.
%
\ep

\begin{Th}
\label{Tcdr}
Let $\bk$ be a bornology on $X$ with closed base. If $X$ has the $\bs$-Hurewicz property then $\cdr(\obs,\obs)$ holds.
\end{Th}
\bp
Let $\{\uc_n:n\in \nb\}$ be a sequence of open $\bs$-covers of $X$. Let $\{Y_n:n\in \nb\}$ be a partition of $\nb$ into infinite subsets. Since  $X$ has the $\bs$-Hurewicz property, ONE has no winning strategy in the $\bs$-Hurewicz game. We define a strategy $\sigma$ for ONE in the $\bs$-Hurewicz game as follows.

Let the first move of ONE be $\sigma(\emptyset)=\uc_1$. TWO responds by choosing a finite set $\vc_1\subseteq \uc_1$. Choose a $k$ for which $1\in Y_k$ and define $\sigma(\vc_1)=\uc_k\setminus \vc_1$. TWO responds by choosing a finite set $\vc_2\subseteq \uc_k$. Suppose that $\vc_n$ has been chosen. Now define $\sigma(\vc_1,\vc_2,\dots,\vc_n)=\uc_m\setminus \{\vc_1\cup \vc_2\cup\dots\cup \vc_n\}$ whenever $n\in Y_m$. TWO responds by choosing a finite set $\vc_{n+1}\subseteq \uc_m$. This define a strategy $\sigma$ for ONE in the $\bs$-Hurewicz game. Since $\sigma$ is not a winning strategy for ONE, consider a play
$$\sigma(\emptyset), \vc_1,\sigma(\vc_1),\vc_2,\dots,\sigma(\vc_1,\vc_2,\dots,\vc_n), \vc_{n+1},\dots$$ which is lost by ONE. Thus for any $B\in \bk$ there exist a $n_0$ and a sequence $\{\delta_n:n\geq n_0\}$ of positive real numbers satisfying $B^{\delta_n}\subseteq U$ for some $U\in \vc_n$ for all $n\geq n_0$. Observe that the collection $\{\vc_n:n\in \nb\}$ is pairwise disjoint by the construction of the strategy. Now for each $m\in \nb$, let $\rc_m=\cup_{n\in Y_m} \vc_{n+1}$. Clearly $\{\rc_n:n\in \nb\}$ is a collection of pairwise disjoint open $\bs$-covers of $X$ witnessing $\cdr(\obs,\obs)$.
\ep

\begin{Th}
\label{Tp4}
Let $\bk$ be a bornology on $X$ with closed base. The  following statements are equivalent.

$(1)$ $X$ has the $\bs$-Hurewicz property.

$(2)$ $X$ satisfies $\Sf(\obs,\obs^{gp})$.

$(3)$ ONE does not have a winning strategy in $\Gf(\obs,\obs^{gp})$.

$(4)$ $X$ satisfies $\obs\rightarrow \ \lceil \obs^{gp}\rceil^2_k$ for each $k\in \nb$, provided $X$ is $\bs$-Lindel\"{o}f.
\end{Th}
\bp
The equivalences of $(1),(2)$ and $(3)$ are due to \cite[Theorem 4.2]{dcpdsd}. We will prove the implications $(3)\Rightarrow (4)$ and $(4)\Rightarrow (2)$.

$(3)\Rightarrow (4)$. This proof is inspired by the arguments used in \cite[Theorem 3]{coocvii}. Let $\uc$ be an open $\bs$-cover of $X$ and let $f:[\uc]^2\rightarrow \{1,\dots,k\}$ be a coloring. Enumerate $\uc$ bijectively as $\{U_n:n\in \nb\}$. Let $\tc_i=\{U_j:j>1,f(\{U_1,U_j\})=i\}$ for $i\in \{1,\dots,k\}$. We can write $\uc\setminus \{U_1\}=\displaystyle{\cup_{i=1}^k} \tc_i$. By Lemma \ref{Lp1}, there is a $i_1\in \{1,\dots,k\}$ such that $\tc_{i_1}$ is an open $\bs$-cover of $X$. Let $\uc_1=\tc_{i_1}$. Applying the same argument inductively we obtain sequences $\{\uc_n:n\in \nb\}$ and $\{i_n:n\in \nb\}$ such that $\uc_{n+1}\subseteq \uc_n$ and $\uc_{n+1}=\{U_j\in \uc_n:j>n+1,f(\{U_{n+1},U_j\})=i_{n+1}\}$ for all $n\in \nb$. Let $\wc_i=\{U_n:i_n=i\}$ for each $i\in \{1,\dots,k\}$. Now  each $\uc_n$ can be expressed as $\displaystyle{\cup_{i=1}^k} (\uc_n\cap \wc_i)$. Again Lemma \ref{Lp1} implies that there must be a $i_n\in \{1,\dots,k\}$ such that $\uc_n\cap \wc_{i_n}$ is an open $\bs$-cover of $X$. Since $\uc_{n+1}\subseteq \uc_n$ for each $n\in \nb$, we can assume that there is a $i_0\in \{1,\dots,k\}$ such that $i_n=i_0$ for all $n\in \nb$.

Define a strategy $\tau$ for ONE in $\Gf(\obs,\obs^{gp})$ as follows. Let the first move of ONE be $\tau(\emptyset)=\uc_1\cap \wc_{i_0}$. TWO selects a finite set $\vc_1\subseteq \tau(\emptyset)$. Let $n_1=\max\{n\in \nb:U_n\in \vc_1\}$. Let the next move of ONE be $\tau(\vc_1)=\uc_{n_1}\cap \wc_{i_0}$. TWO selects a finite set $\vc_2\subseteq \tau(\vc_1)$. Clearly $\vc_1\cap \vc_2=\emptyset$. Let $n_2=\max\{n\in \nb:U_n\in \vc_2\}$. Again let $\tau(\vc_1,\vc_2)=\uc_{n_2}\cap \wc_{i_0}$. TWO selects a finite set $\vc_3\subseteq \tau(\vc_1,\vc_2)$ and so on. Clearly $\{\vc_r:r\in \nb\}$ are pairwise disjoint.

This defines $\tau$ a strategy for ONE in $\Gf(\obs,\obs^{gp})$. Since $\tau$ is not a winning strategy for ONE, consider a play
\[\tau(\emptyset),\vc_1,\dots,\tau(\vc_1,\dots,\vc_r), \vc_{r+1},\dots.\]
which is lost by ONE. So $\displaystyle{\cup_{r\in \nb}}\vc_r\in \obs^{gp}$. Let $\vc=\cup_{r\in \nb}\vc_r$. Then $\vc\in \obs^{gp}$ and $\{\vc_r:r\in \nb\}$ is a partition of $\vc$ into pairwise disjoint finite sets such that for any $V,W\in \vc$ where $V\in \vc_r$ and $W\in \vc_s$ with $r\neq s$, $f(\{V,W\})=i_0$. This shows that  $(4)$ holds.

$(4)\Rightarrow (2)$ Clearly $X$ satisfies $\obs\rightarrow \ \lceil \obs\rceil^2_2$. By Theorem \ref{Tp1}, $X$ satisfies $\Sf(\obs,\obs)$. Now we show that every countable open $\bs$-cover of $X$ is $\bs$-groupable. Let $\uc$ be a countable open $\bs$-cover of $X$. Let $\{\uc_n:n\in \nb\}$ be a partition of $\uc$ into nonempty finite sets. Define a coloring $f:[\uc]^2\rightarrow \{1,2\}$ by $f(\{U,V\})=1$ if $U,V\in \uc_n$ for some $n\in \nb$ and $f(\{U,V\})=2$ otherwise. By $(4)$, there is a $\vc\in \obs^{gp}$ with $\vc\subseteq \uc$. So there is a sequence $\{\vc_n:n\in \nb\}$ of pairwise disjoint finite subsets of $\vc$ witnessing the $\bs$-groupability of $\vc$. Since $\uc$ is countable, the elements of $\uc\setminus \vc$ can be distributed among $V_n$'s so that $\{V_n:n\in \nb\}$ witnesses the $\bs$-groupability of $\uc$. Hence $X$ satisfies $\Sf(\obs,\obs^{gp})$.

\ep

Our next results concerns with the $\bs$-Hurewicz property of product spaces.
For a metric space $(X,d)$ one can consider the product space $X^n$ endowed with the product metric $d^n$ defined as $$d^n((x_1,\dotsc, x_n),(y_1,\dotsc, y_n))=\max\{d(x_1,y_1),\dotsc, d(x_n,y_n)\}.$$

Let $\bk$ be a bornology on $X$ with closed base $\bk_0$. In \cite{hh}, it has been shown that the collection $\{B^n:B\in \bk\}$ generates a bornology on $X^n$. We denote that bornology on $X^n$ by $\bk^n$. Further it can be easily verified that for any $n\in \nb$ and $\delta>0$, $(B^\delta)^n=(B^n)^\delta$ which we will need repeatedly. We start with the following simple observation followed by a result on selection principles (Theorem \ref{Tpdt1}) which though not related to the main topic of this section is interesting in its own right.

\begin{Lemma}
\label{Lpdt}
Let $\bk$ be a bornology on $X$ with a compact base $\bk_0$. If $\uc$ is an open $\bns$-cover of $X^n$, then there exists an open $\bs$-cover $\vc$ of $X$ such that $\{V^n:V\in \vc\}$ is an open $\bns$-cover of $X^n$ which refines $\uc$.
\end{Lemma}
\bp
Let $\uc$ be an open $\bns$-cover of $X^n$. Let $B\in \bk_0$. Then for $B^n\in \bn$ there exist a $\delta>0$ and $U\in \uc$ such $(B^n)^\delta\subseteq U$. Since $B$ is a compact subset of $X$ and $U$ is open in $X^n$ containing $B^n$, we use Wallace Theorem \cite{engelking} to find an open set $V_B$ in $X$ such that $B^n\subseteq V_B^n\subseteq U$. Consider the collection $\vc=\{V_B:B\in \bk_0\}$. We intend to show that $\vc$ is an open $\bs$-cover of $X$. For $B\in \bk_0$, we have $B\subseteq V_B$. Since $V_B$ is open and $B$ is compact with $B\subseteq V_B$, there is a $\delta>0$ such that $B^\delta\subseteq V_B$.
%for each $x\in B$ there exists a $\delta_x>0$ such that $S(x,\delta_x)\subseteq V_B$. Consider the collection $\{S(x,\frac{\delta_x}{2}):x\in B\}$. As $B$ is compact, there are $x_1,\dotsc,x_k$ such that $B\subseteq \displaystyle{\cup_{i=1}^k} S(x_i,\frac{\delta_{x_i}}{2})\subseteq \displaystyle{\cup_{i=1}^k} S(x_i,\delta_{x_i})\subseteq V_B$. Choose $\delta=\min\{\frac{\delta_{x_i}}{2}:i=1,\dotsc,k\}$. It is easy to observe that $B^\delta\subseteq V_B$.
%For $x\in B^\delta$, there exists a $b\in B$ such that $d(x,b)<\delta$. Again for $b\in B$, $d(b,x_i)<\frac{\delta_{x_i}}{2}$ for some $i\in \{1,\dotsc,k\}$. So $d(x,x_i)<\delta+\frac{\delta_{x_i}}{2}<\delta_{x_i}$ i.e, $x\in S(x_i,\delta_{x_i})$ for some $i\in \{1,\dotsc,k\}$. So $B^\delta\subseteq \cup_{i=1}^k S(x_i,\delta_{x_i})\subseteq V_B$.
Therefore $\vc=\{V_B:B\in \bk_0\}$ is an open $\bs$-cover of $X$.

Again for $B\in \bk_0$, $(B^\delta)^n\subseteq V_B^n$ i.e. $(B^n)^\delta\subseteq V_B^n$ and $V_B\subseteq U$ for some $U\in \uc$ implying that $\{V^n:V\in \vc\}$ is an open $\bns$-cover of $X^n$ as well as refines $\uc$.
\ep

\begin{Th}
\label{Tpdt1}
Let $\Pi\in \{\S1,\Sf\}$ and $\pc,\qc\in \{\oc,\Gamma\}$. Let $\bk$ be a bornology on $X$ with compact base. The following statements are equivalent.

$(1)$ $X$ satisfies $\Pi(\pc_{\bs},\qc_{\bs})$.

$(2)$ $X^n$ satisfies $\Pi(\pc_{\bns},\qc_{\bns})$ for each $n\in \nb$.
\end{Th}
\bp
We only present the proof for the case $\Pi=\S1$, $\pc=\oc$ and $\qc=\oc$ as other cases follow analogously.

$(1)\Rightarrow (2)$. Let $\{\uc_m:m\in \nb\}$ be a sequence of open $\bns$-covers of $X^n$. By Lemma \ref{Lpdt}, for each $\uc_m$ there exists an open $\bs$-cover $\vc_m$ of $X$ such that $\{V^n:V\in \vc_m\}$ refines $\uc_m$. Apply $\S1(\obs,\obs)$ to the sequence $\{\vc_m:m\in \nb\}$ to choose $V_m\in \vc_m$ for each $m\in \nb$ so that  $\{V_m:m\in \nb\}$ becomes an open $\bs$-cover of $X$. Now for each $V_m$ we can choose a $U_m\in \uc_m$ such that $V_m^n\subseteq U_m$. We will show that $\{U_m:m\in \nb\}$ is an open $\bns$-cover of $X^n$. Let $B^n\in \bn$. For $B\in \bk$ there exist a $\delta>0$ and a $V_m$ such that $B^\delta\subseteq V_m$. So $(B^\delta)^n\subseteq V^n_m$ i.e. $(B^n)^\delta\subseteq V^n_m\subseteq U_m$. Therefore $\{U_m:m\in \nb\}$ is an open $\bns$-cover of $X^n$. Hence $X^n$ satisfies $\S1(\obns,\obns)$.

$(2)\Rightarrow (1)$. Let $\{\uc_m:m\in \nb\}$ be a sequence of open $\bs$-covers of $X$. For each $m\in \nb$, let $\uc_m'=\{U^n:U\in \uc_m\}$. Clearly $\uc_m'$'s are open $\bns$-cover of $X^n$. Now consider the sequence $\{\uc_m':m\in \nb\}$ of open $\bns$-covers of $X^n$. Since $X^n$ satisfies $\S1(\obns,\obns)$, there is an open $\bns$-cover $\{U^n_m:m\in \nb\}$ of $X^n$ with $U^n_m\in \uc_m'$ for each $m\in \nb$. We will show that $\{U_m:m\in \nb\}$ with $U_m\in \uc_m$ for $m\in \nb$ is an open $\bs$-cover of $X$. Let $B\in \bk$. Then for $B^n\in \bn$ there exist a $\delta>0$ and $U^n_m$ such that $(B^n)^\delta\subseteq U^n_m$ i.e. $(B^\delta)^n\subseteq U^n_m$ i.e. $B^\delta\subseteq U_m$. So $\{U_m:m\in \nb\}$ is an open $\bs$-cover of $X$. Hence $X$ satisfies $\S1(\obs,\obs)$.
\ep

\begin{Th}
\label{Tpdt2}
Let $\bk$ be a bornology on $X$ with compact base. The following statements are equivalent.

$(1)$ $X$ has the $\bs$-Hurewicz property.

$(2)$ $X^n$ has the $\bns$-Hurewicz property for each $n\in \nb$.
\end{Th}
\bp
We only present proof of $(1)\Rightarrow (2)$.\\
$(1)\Rightarrow (2)$. Let $\{\uc_k:k\in \nb\}$ be a sequence of open $\bns$-covers of $X^n$. By Lemma \ref{Lpdt}, for each $k$ there exists an open $\bs$-cover $\vc_k$ of $X$ such that $\{V^n:V\in \vc_k\}$ refines $\uc_k$. Consider the sequence $\{\vc_k:k\in \nb\}$. By $(1)$, there is a sequence $\{\wc_k:k\in \nb\}$ of finite sets with $\wc_k\subseteq \vc_k$ for each $k$ such that for $B\in \bk$ there exist a $m_0\in \nb$ and a sequence $\{\delta_m:m\ge  m_0\}$ of positive reals satisfying $B^{\delta_m}\subseteq V$ for some $V\in \wc_k$ for all $m\ge  m_0$. Now for each $k$ we can choose a finite subset $\zc_k$ of $\uc_k$ such that for each $V\in \wc_k$ there is a $U\in \zc_k$ with $V^n\subseteq U$. We will show that $\{\zc_k:k\in \nb\}$ witnesses the $\bns$-Hurwicz property of $X^n$. Let $B^n\in \bn$. Note that for $B\in \bk$ there already exist a $p_0\in \nb$ and a sequence $\{\delta_p:p\ge  p_0\}$ of positive reals satisfying $B^{\delta_p}\subseteq V$ for some $V\in \wc_k$ for all $p\ge  p_0$ i.e. $(B^n)^{\delta_p}\subseteq U$ for some $U\in \zc_k$ for all $p\ge  p_0$. Hence $X^n$ has the $\bns$-Hurewicz property.

%$(2)\Rightarrow (1)$. Let $\{\uc_k:k\in \nb\}$ be a sequence of open $\bs$-covers of $X$. For each $k\in \nb$, let $\uc^n_k=\{U^n:U\in \uc_k\}$. Clearly $\uc^n_k$'s are open $\bns$-cover of $X^n$ for each $k\in \nb$. Consider the sequence $\{\uc^n_k:k\in \nb\}$. By $(2)$, there is a sequence $\{\vc_k:k\in \nb\}$ of finite sets with $\vc_k\subseteq \uc^n_k$ for each $k\in \nb$ such that for $B^n\in \bn$ there are a $m_0\in \nb$ and a sequence $\{\delta_m:m \geq m_0\}$ of positive reals satisfying $(B^n)^{\delta_m}\subseteq U^n$ for some $U^n\in \vc_k$ for all $m\ge  m_0$. Let $\wc_k=\{U\in \uc_k:U^n\in \vc_k\}$ for each $k\in \nb$. Now $(B^{\delta_m})^n\subseteq U^n$ for some $U^n\in \vc_k$ for all $m\ge  m_0$ i.e. $B^{\delta_m}\subseteq U$ for some $U\in \wc_k$ for all $m\ge  m_0$. This shows that $\{\wc_k:k\in \nb\}$ witnesses the $\bs$-Hurewicz property of $X$. Hence $(1)$ holds.
\ep

Combining Theorem \ref{Tp4} and Theorem \ref{Tpdt2}, we obtain the following.
\begin{Th}
\label{Tpdt3}
Let $\bk$ be a bornology on $X$ with compact base. The  following statements are equivalent.

$(1)$ $X^n$ has the $\bns$-Hurewicz property for each $n\in \nb$.

$(2)$ $X$ satisfies $\Sf(\obs,\obs^{gp})$.

$(3)$ ONE does not have a winning strategy in $\Gf(\obs,\obs^{gp})$ on $X$.

$(4)$ $X$ satisfies $\obs\rightarrow \ \lceil \obs^{gp}\rceil^2_k$ for each $k\in \nb$, provided $X$ is $\bs$-Lindel\"{o}f.
\end{Th}

\subsection{The strong $\bk$-Gerlits-Nagy property and some observations}
It is well known that the Gerlits-Nagy property was introduced in the seminal paper of the authors \cite{gn} as a property stronger than the classical Hurewicz property and had been extensively investigated since then (see \cite{sur3,tz} for example). In particular in \cite{coocvii} certain new characterizations of this property were established using the notion of groupability of open covers. In this section we introduce the notion of strong $\bk$-Gerlits-Nagy property (which has not been investigated at all in bornological settings) and while defining the concept, follow the line of \cite{coocvii} which seems more effective for our purpose.
\begin{Def}
\label{Dgn}
Let $\bk$ be a bornology on $X$ with closed base. $X$ is said to have the strong $\bk$-Gerlits-Nagy property (in short, $\bs$-Gerlits-Nagy property) if $X$ satisfies the selection principle $\S1(\obs,\obs^{gp})$.
\end{Def}
It is clear that  $\S1(\obs,\Gbs)$ implies $\S1(\obs,\obs^{gp})$ which shows that every $\gbs$-set has the $\bs$-Gerlits-Nagy property. Again $\S1(\obs,\obs^{gp})$ evidently implies $\Sf(\obs,\obs^{gp})$ which in turn assures the $\bs$-Hurewicz property by \cite[Theorem 4.2]{dcpdsd}. Therefore as in the classical case, if $X$ has the $\bs$-Gerlits-Nagy property then $X$ has the $\bs$-Hurewicz property and further it satisfies $\S1(\obs,\obs)$. Below we prove certain results in line of \cite{coocvii} where the proofs are done with suitable modifications as is necessary for bornological structures.

The following example shows that the real line associated with a bornology has the $\bs$-Gerlits-Nagy property.
\begin{Ex}
\label{H1}
Consider the real line $X=\rb$ with the Euclidean metric $d$ and the bornology $\bk$ generated by $\{(-x,x):x>0\}$. We show that $X$ has the $\bs$-Gerlits-Nagy property. To see this, let $\{\uc_n:n\in \nb\}$ be a sequence of open $\bs$-covers of $X$. Then clearly for each $k\in\nb$, there is a $U\in\uc_n$ such that $(-k,k)\subseteq U$ (by choosing $(-k,k)\in \bk$).
%For $B=(-k,k)\in \bk$, where $k$ is a positive integer, there exist a $\delta>0$ and $U\in \uc_n$ such that $B^\delta\subseteq U$ so $(-k,k)\subseteq U$.

Consider a sequence of positive integers $k_1< k_2< \cdots$. For each $n\in \nb$, choose a $U_n\in \uc_n$ such that $(-k_n,k_n)\subseteq U_n$. We show that $\{U_n:n\in \nb\}\in \obs^{gp}$. For this we choose $\vc_n=\{U_n\}$ for each $n\in \nb$. Now we show that $\{\vc_n:n\in\nb\}$ is the sequence of pairwise disjoint finite sets witnessing the $\bs$-groupability of $\{U_n:n\in \nb\}$.
Let $B\in\bk$. Clearly $\uc=\{(-n,n):n\in \nb\}$ is a $\gamma_{\bs}$-cover of $X$. Consequently we can find a $n_0\in \nb$ and a sequence of positive real numbers $\{\delta_n:n\geq  n_0\}$ such that $B^{\delta_n}\subseteq (-n,n)$ for all $n\geq  n_0$ i.e. $B^{\delta_n}\subseteq (-k_n,k_n)\subseteq U_n$ for $U_n\in \vc_n$ for all $n\geq  n_0$. Therefore $\{U_n:n\in \nb\}\in \obs^{gp}$ and so $X$ satisfies $\S1(\obs,\obs^{gp})$. Hence $X$ has the $\bs$-Gerlits-Nagy property.
\end{Ex}

\begin{Th}
\label{Tgn-1}
Let $\bk$ be a bornology on $X$ with closed base. The following statements are equivalent.\\
\noindent$(1)$ $X$ has the $\bs$-Gerlits-Nagy property.\\
\noindent$(2)$ $X$ has the $\bs$-Hurewicz property as well as it satisfies $\S1(\obs,\obs)$.
\end{Th}
\bp
\noindent$(1)\Rightarrow (2)$. By $(1)$, $X$ satisfies $\S1(\obs,\obs^{gp})$ which implies that $X$ satisfies $\Sf(\obs,\obs^{gp})$ as well as $\S1(\obs,\obs)$. Again by \cite[Theorem 4.2]{dcpdsd}, $X$ has the $\bs$-Hurewicz property. Hence $(2)$ holds.

\noindent$(2)\Rightarrow (1)$. Let $\{\uc_n:n\in \nb\}$ be a sequence of open $\bs$-covers
of $X$. Apply $\S1(\obs,\obs)$ to $\{\uc_n:n\in \nb\}$ to choose a $U_n\in \uc_n$ for each
$n\in \nb$ such that $\{U_n:n\in \nb\}$ is an open $\bs$-cover of $X$. Using that $X$ has the $\bs$-Hurewicz property and \cite[Lemma 4.2]{dcpdsd}, $\{U_n:n\in \nb\}$ is a $\bs$-groupable cover of $X$. Thus $X$ satisfies $\S1(\obs,\obs^{gp})$.
\ep

\begin{Prop}
\label{Pgn}
Let $\bk$ be a bornology on $X$ with a closed base $\bk_0$. If $|\bk_0|<\add(\mc)$ then $X$ has the $\bs$-Gerlits-Nagy property.
\end{Prop}
\bp
Since $\add(\mc)=\min\{\bb,\cov(\mc)\}$ and $|\bk_0|<\add(\mc)$, so $|\bk_0|<\cov(\mc)$ as well as $|\bk_0|<\bb$. By Theorem \ref{Tcard4}, $|\bk_0|<\cov(\mc)$ implies that $X$ satisfies $\S1(\obs,\obs)$. Again in view of $|\bk_0|<\bb$ and Theorem \ref{Thb}, we can conclude  that $X$ has the $\bs$-Hurewicz property. Hence $X$ has the $\bs$-Gerlits-Nagy property.
\ep

\begin{Th}
\label{Tgn-2}
Let $\bk$ be bornology on $X$ with closed base. The following statements are equivalent.

$(1)$ $X$ has the $\bs$-Gerlits-Nagy property.

$(2)$ $X$ satisfies $\obs\rightarrow (\obs^{gp})^n_k$ for each $n,k\in \nb$.
\end{Th}
\bp
$(1)\Rightarrow (2)$. Let $\uc$ be an open $\bs$-cover of $X$ and  $f:[\uc]^n\rightarrow \{1,\dots,k\}$ be a coloring.  Since $X$ satisfies $\S1(\obs,\obs)$, $X$ also satisfies $\obs\rightarrow (\obs)^n_k$ by Theorem \ref{Tp2}. Consequently there are a $\vc\subseteq \uc$ with $\vc\in \obs$ and a $i\in \{1,\dots,k\}$ such that for each $V\in [\vc]^n$, $f(V)=i$. By $(1)$, we can find a countable open $\bs$-cover $\vc'\subseteq \vc$ which is $\bs$-groupable. Thus we have a $\vc'\in \obs^{gp}$ and a $i\in \{1,\dots,k\}$ such that for each $V\in [\vc']^n$, $f(V)=i$ holds.

$(2)\Rightarrow (1)$ Clearly $X$ satisfies $\obs\rightarrow (\obs)^2_2$. By Theorem \ref{Tp2}, $X$ satisfies $\S1(\obs,\obs)$. Now we show that every countable open $\bs$-cover of $X$ is $\bs$-groupable. Let $\uc$ be a countable open $\bs$-cover of $X$. Let $\{\uc_n:n\in \nb\}$ be a partition of $\uc$ into nonempty finite sets. Define a coloring $f:[\uc]^2\rightarrow \{1,2\}$ by $f(\{U,V\})=1$ if $U,V\in \uc_n$ for some $n\in \nb$ and $f(\{U,V\})=2$ otherwise. By $(2)$, there is a $\vc\in \obs^{gp}$ with $\vc\subseteq \uc$. Therefore $\uc$ is $\bs$-groupable. Hence $X$ satisfies $\S1(\obs,\obs^{gp})$.

\ep

\begin{Th}
\label{Tgn-3}
Let $\bk$ be bornology on $X$ with compact base. The following statements are equivalent.

$(1)$ $X$ has the $\bs$-Gerlits-Nagy property.

$(2)$ $X^n$ has the $(\bk^n)^s$-Gerlits-Nagy property for each $n\in \nb$.
\end{Th}
\bp
We only prove $(1)\Rightarrow (2)$.

Let $\{\uc_k:k\in \nb\}$ be a sequence of open $(\bk^n)^s$-covers of $X^n$. By Lemma \ref{Lpdt}, for each $\uc_k$ there exists an open $\bs$-cover $\vc_k$ of $X$ such that $\{V^n:V\in \vc_k\}$ refines $\uc_k$. Apply  $\S1(\obs,\obs^{gp})$ to $\{\vc_k:k\in \nb\}$ to choose a $V_k\in \vc_k$ for each $k\in \nb$ such that $\{V_k:k\in \nb\}$ is a $\bs$-groupable cover of $X$. Now for each $V_k$ there is a $U_k\in \uc_k$ such that $V^n_k\subseteq U_k$. We  will show that $\{U_k:k\in \nb\}$ is a $(\bk^n)^s$-groupable cover of $X^n$.

 As $\{V_k:k\in \nb\}\in \obs^{gp}$, there is a sequence $\{\wc_p:p\in \nb\}$ of pairwise disjoint finite sets such that for $B\in \bk$ there exist a $p_0\in \nb$ and a sequence $\{\delta_p:p\geq p_0\}$ of positive reals satisfying $B^{\delta_p}\subseteq V_k$ for some $V_k\in \wc_p$ for all $p\geq p_0$. Let $\zc_p$ for $p\in \nb$ be a finite subset of $\{U_k:k\in \nb\}$ such that for each $V_k\in \wc_p$, $V^n_k\subseteq U_k$ for some $U_k\in \zc_p$. Now as before we will have  $(B^n)^{\delta_p} = (B^{\delta_p})^n\subseteq V^n_k\subseteq U_k$  for some $U_k\in \zc_p$ for all $p\geq p_0$. Therefore $\{U_k:k\in \nb\}$ is a $(\bk^n)^s$-groupable cover of $X^n$ witnessing $\S1(\obns,\obns^{gp})$ implying the $(\bk^n)^s$-Gerlits-Nagy property  of $X^n$.
\ep

\subsection{ Observation on product spaces $X\times Y$}
We end Section 4 with two more observations regarding both the properties in product spaces. Let $\bk$ be a bornology on $X$ and $(Y,\rho)$ be another metric space. There is a natural bornology $\widehat{\bk}$ on $(X\times Y,d\times \rho)$ induced by $\bk$ defined as $\widehat{\bk}=\{C\subseteq X\times Y:\pi_X(C)\in \bk\}$ \cite{bl}, where $\pi_X:X\times Y\rightarrow X$ is the projection map. A base for $\widehat{\bk}$ is $\{B\times Y:B\in \bk\}$.

\begin{Lemma}
\label{Lpdt2}
Let $\bk$ be a bornology on $X$ with compact base and $(Y,\rho)$ be another compact metric space. If $\uc$ is an open $(\widehat{\bk})^s$-cover of $X\times Y$, then there exists an open $\bs$-cover $\vc$ of $X$ such that $\{V\times Y:V\in \vc\}$ refines $\uc$.
\end{Lemma}
The proof is analogous to the proof of Lemma \ref{Lpdt}.

\begin{Th}
\label{Tpdt4}
Let $\bk$ be a bornology on $X$ with compact base and $(Y,\rho)$ be another compact metric space. The following statements hold.

$(1)$ $X$ satisfies $\S1(\obs,\Gbs)$ if and only if $X\times Y$ satisfies $\S1(\obhs,\Gbhs)$.

$(2)$ $X$ has the $\bs$-Hurewicz property if and only if $X\times Y$ has the $(\widehat{\bk})^s$-Hurewicz property.

$(3)$ $X$ has the $\bs$-Gerlits-Nagy property if and only if $X\times Y$ has the $(\widehat{\bk})^s$-Gerlits-Nagy property.
\end{Th}

Let $(Y,d_Y)$ be a subspace of $(X,d)$ where $d_Y$ is the induced metric on $Y\subseteq X$. Let $\bk$ be a bornology on $X$. It can be easily checked that $\bk_Y=\{B\cap Y:B\in \bk\}$ is a bornology on $Y$. If $Y$ is closed and $\bk$ has a compact base, then $\bk_Y$ has a compact base. The following result needs no further explanations.

\begin{Lemma}
\label{TY}
Let $\bk$ be a bornology on $X$ with compact base and $Y$ be a closed subset of $X$. The following statements are true.

$(1)$ If $X$ satisfies $\S1(\obs,\Gbs)$, then $Y$ satisfies $\S1(\oc_{\bk^s_Y},\Gamma_{\bk^s_Y})$.

$(2)$ If $X$ has the $\bs$-Hurewicz property, then $Y$ has the $\bk^s_Y$-Hurewicz property.

$(3)$ If $X$ has the $\bs$-Gerlits-Nagy property, then $Y$ has the $\bk^s_Y$-Gerlits-Nagy property.
\end{Lemma}

%%%%%%%%%%%%%%%%%%%%%%%%%%%%%added%%%%%%%%%%%%%%%%%
%%%%%%%%%%%%%%%%%%%%%%%%%%%%%%%%%%%%%%%%%%%%%%%%%%%%%%

\section{Results in function spaces}

\subsection{Further observations regarding selection principles and $C(X)$}
In this section we continue our investigation of selection principles in function spaces which we had started in \cite{dcpdsd} and establish  some new results which have not been discussed in earlier articles.

As before, let $\bk$ be a bornology on $(X,d)$ with closed base and $(Y,\rho)$ be another metric space. For $f\in C(X,Y)$, the neighbourhood of $f$ with respect to the topology $\tau_{\bk}^s$ of strong uniform convergence is denoted by
\[ [B,\varepsilon ]^s(f)=\{g\in C(X,Y):\exists \delta>0,\rho(f(x),g(x))<\varepsilon , \forall x\in B^\delta\}, \]
for $B\in\bk,\,\varepsilon >0$ (\!\cite{bl,cmh}).

The symbol $\underline{0}$ denotes the  zero function in $(C(X), \tau_{\bk}^s)$. The space $(C(X), \tau_{\bk}^s)$ is homogeneous and so it is enough to concentrate at the point $\underline{0}$ when dealing with local properties of this function space.

%%%%%%%%%%%%%%%%%%%%%%%%%%%%%%%added%%%%%%%%%%%%%%%%%%%
%%%%%%%%%%%%%%%%%%%%%%%%%%%%%%%%%%%%%%%%%%%%%%%%%%%%%%%
%%%%%%%%%%%%%%%%%%%%%%%%%%%%%%%%%%%%%%%%%%%%%%%%%%%%%
We start with recollecting the following result from \cite{cmk} which will be useful in our context.
\begin{Lemma}{(\cite[Lemma 2.2]{cmk})}
\label{koc}
Let  $\bk$ be a bornology on the metric space $(X,d)$ with closed base. The following statements hold.

\noindent $(a)$ Let $\uc$ be an open $\bs$-cover of $X$. If $A = \{f \in  C(X): \exists U \in\uc, f(x) = 1 \,\,\text{for all}\,\, x \in X \setminus U\}$. Then $\underline{0} \in \overline{A} \setminus A$ in $(C(X), \tau_\bk^s)$.\\
\noindent $(b)$ Let $A \subseteq (C(X), \tau_\bk^s)$ and let $\uc = \{f^{-1}(-\frac{1}{n},\frac{1}{n}): f \in A\}$, where $n\in\nb$. If $\underline{0} \in \overline{A}$ and $X\notin \uc$, then $\uc$ is an  open $\bs$-cover of $X$.
\end{Lemma}

The following exemplary observation about $\gbs$-covers is also useful.
%%%%%
%%change the statement of lemma in terms of epsilon. Does the converse hold?
%%%%%
\begin{Lemma}
\label{L8}
Let $\bk$ be a bornology on $X$ with closed base. Let $A=\{f_n:n\in \nb\}$ be a sequence of functions in $(C(X),\tau^s_\bk)$. If $A\in \Sigma_{\underline{0}}$ then for any neighbourhood $U$ of $0$ in the real line, $\{f^{-1}_n(U):n\in \nb\}$ is a $\gbs$-cover of $X$.
\end{Lemma}
\bp
As $U$ is a neighbourhood of $0$, we can find a $\varepsilon>0$ such that $(-\varepsilon,\varepsilon)\subseteq U$. Let $B\in \bk$. Consider the neighbourhood $[B,\varepsilon]^s(\underline{0})$ of $\underline{0}$. Since $\{f_n:n\in \nb\}$ converges to $\underline{0}$ with respect to $\tau^s_\bk$, we can choose a $n_0\in \nb$ such that $f_n\in [B,\varepsilon]^s(\underline{0})$ for all $n\geq  n_0$. This means that for each $n\geq  n_0$ there exists a $\delta_n>0$ such $B^{\delta_n}\subseteq f^{-1}_n(-\varepsilon,\varepsilon)\subseteq f^{-1}_n(U)$ i.e. for the sequence $\{\delta_n:n\geq  n_0\}$ of positive real numbers we have $B^{\delta_n}\subseteq f^{-1}_n(U)$ for all $n\geq  n_0$. Hence $\{f^{-1}_n(U):n\in \nb\}$ is a $\gbs$-cover of $X$.
\ep

%The next result is useful in what follows.
%\begin{Lemma}
%\label{L6}
%Let $\bk$ be a bornology on $X$ with closed base. If $\{f_n:n\in \nb\}$ and $\{g_n:n\in \nb\}$ are two sequence of functions in $(C(X),\tau^s_\bk)$ converging to $\underline{0}$ with respect to $\tau^s_\bk$, then $\{f_n+ g_n:n\in \nb\}$ also converges to $\underline{0}$ with respect to $\tau^s_\bk$.
%\end{Lemma}

\begin{Prop}
\label{Tc}
Let $\bk$ be a bornology on $X$ with closed base. If $X$ satisfies $\S1(\Gbs, \Gbs)$ then $(C(X),\tau^s_\bk)$ satisfies $\S1(\Sigma_{\underline{0}},\Sigma_{\underline{0}})$.
\end{Prop}
\bp
Let $\{A_n:n\in \nb\}$ be a sequence of elements in $\Sigma_{\underline{0}}$, where $A_n=\{f_{n,k}:k\in \nb\}$ for $n\in \nb$. Consider the set $\uc_n=\{f^{-1}_{n,k}(-\frac{1}{n},\frac{1}{n}):k\in \nb\}$. By Lemma \ref{L8}, each $\uc_n$ is a $\gbs$-cover of $X$. Now apply $\S1(\Gbs, \Gbs)$ to $\{\uc_n:n\in \nb\}$ to choose a $f^{-1}_{n,k_n}(-\frac{1}{n},\frac{1}{n})\in \uc_n$ for each $n\in \nb$ such that $\{f^{-1}_{n,k_n}(-\frac{1}{n},\frac{1}{n}):n\in \nb\}\in \Gbs$. We claim that $\{f_{n,k_n}:n\in \nb\}$ converges to $\underline{0}$. Let $[B,\varepsilon]^s(\underline{0})$ be neighbourhood of $\underline{0}$ where $B\in \bk$ and $\varepsilon>0$. First choose $n_1\in \nb$ such that $\frac{1}{n_1}<\varepsilon$. Since $\{f^{-1}_{n,k_n}(-\frac{1}{n},\frac{1}{n}):n\in \nb\}$ is a $\gbs$-cover, there exist a $n_0$ and a sequence $\{\delta_n:n\geq n_0\}$ of positive real numbers such that $B^{\delta_n}\subseteq f^{-1}_{n,k_n}(-\frac{1}{n},\frac{1}{n})$ for all $n\geq n_0$. Clearly $B^{\delta_n}\subseteq f^{-1}_{n,k_n}(-\varepsilon,\varepsilon)$ for all $n\geq n_2$, where $n_2=\max\{n_0,n_1\}$. This shows that $\{f_{n,k_n}:n\in \nb\}$ converges to $\underline{0}$ with respect to $\tau^s_\bk$. Therefore $(C(X),\tau^s_\bk)$ satisfies $\S1(\Sigma_{\underline{0}},\Sigma_{\underline{0}})$.
\ep

\begin{Th}
\label{Tf1}
Let $\bk$ be a bornology on $X$ with closed base. The following statements are equivalent.\\
\noindent$(1)$ ONE has no winning strategy in $\G1(\obs,\Gbs)$ on $X$. \\
\noindent$(2)$ ONE has no winning strategy in $\G1(\Omega_{\underline{0}},\Sigma_{\underline{0}})$ on $(C(X),\tau^s_\bk)$.
\end{Th}
\bp
Let $\sigma$ be a strategy for ONE in $\G1(\Omega_{\underline{0}},\Sigma_{\underline{0}})$ on $(C(X),\tau^s_\bk)$. We now use $\sigma$ to define a strategy $\psi$ for ONE in $\G1(\obs,\Gbs)$ on $X$ as follows.

The first move of ONE in  $\G1(\Omega_{\underline{0}},\Sigma_{\underline{0}})$ is $\sigma(\emptyset)$. Let $\uc_1=\{f^{-1}(-1,1):f\in \sigma(\emptyset)\}$ and assume that $X\not\in \uc_1$. By Lemma \ref{koc}, $\uc_1$ is an open $\bs$-cover of $X$. Define $\psi(\emptyset)=\uc_1$, the first move of ONE in $\G1(\obs,\Gbs)$. TWO responds with $U_1=f^{-1}_1(-1,1)$. Let the move of TWO in $\G1(\Omega_{\underline{0}},\Sigma_{\underline{0}})$ be $f_1$. To define $\psi(U_1)$, we look at the move $\sigma(f_1)$ of ONE in  $\G1(\Omega_{\underline{0}},\Sigma_{\underline{0}})$. Consider $\uc_2=\{f^{-1}(-\frac{1}{2},\frac{1}{2}):f\in \sigma(f_1)\}$, which is again an open $\bs$-cover of $X$ by Lemma \ref{koc}. Define $\psi(U_1)=\uc_2$. TWO responds with $U_2=f^{-1}_2(-\frac{1}{2},\frac{1}{2})$. Let TWO's move in $\G1(\Omega_{\underline{0}},\Sigma_{\underline{0}})$ be $f_2$ and so on.

This defines a strategy $\psi$ for ONE in $\G1(\obs,\Gbs)$. Since $\psi$ is not a winning strategy, consider a $\psi$-play
\[\psi(\emptyset), U_1, \psi(U_1), U_2, \dots\]
 which is lost by ONE in $\G1(\obs,\Gbs)$. Thus $\{U_n:n\in \nb\}$ is a $\gbs$-cover of $X$, where $U_n=f^{-1}_n(-\frac{1}{n},\frac{1}{n})$  for each $n\in \nb$. Consequently $\{f_n:n\in \nb\}$ is converges to $\underline{0}$ with respect to $\tau^s_\bk$.

Now correspond to the $\psi$-play there is also a $\sigma$-play
\[\sigma(\emptyset), f_1, \sigma(f_1), f_2,\dotsc\]
 in $\G1(\Omega_{\underline{0}},\Sigma_{\underline{0}})$ and $\{f_n:n\in \nb\}\in \Sigma_{\underline{0}}$. Hence $\sigma$ is not a winning strategy for ONE in $\G1(\Omega_{\underline{0}},\Sigma_{\underline{0}})$.\\

$(2)\Rightarrow (1)$. Let $\psi$ be a strategy for ONE in $\G1(\obs,\Gbs)$ on $X$. We define a strategy $\sigma$ for ONE in $\G1(\Omega_{\underline{0}},\Sigma_{\underline{0}})$ on $(C(X),\tau^s_\bk)$ as follows.

The first move of ONE in $\G1(\obs,\Gbs)$ is $\psi(\emptyset)=\uc_1$ (say). Since $\uc_1\in \obs$, for $B\in \bk$ there are a $\delta>0$ and a $U\in \uc_1$ such that $B^{2\delta}\subseteq U$. Let $\uc_{1,B}=\{U\in \uc_1:B^{2\delta}\subseteq U\}$. For each $U\in \uc_{1,B}$ choose a $f_{B,U}\in C(X)$ such that $f_{B,U}(B^\delta)=\{0\}$ and $f_{B,U}(X\setminus U)=\{1\}$. Consider the collection $A_1=\{f_{B,U}:B\in \bk,U\in \uc_{1,B}\}$. Clearly $A_1\in \Omega_{\underline{0}}$. Now define $\sigma(\emptyset)=A_1$, the first move of ONE in $\G1(\Omega_{\underline{0}},\Sigma_{\underline{0}})$. TWO responds by choosing $f_{B_1,U_1}\in A_1$. Let $U_1$ be the TWO's response in $\G1(\obs,\Gbs)$. Now let $\psi(U_1)=\uc_2$. We similarly construct $A_2=\{f_{B,U}:B\in \bk,U\in \uc_{2,B}\}$ which is in $\Omega_{\underline{0}}$. Define $\sigma(f_{B_1,U_1})=A_2$. TWO responds by choosing $f_{B_2,U_2}\in A_2$. Let $U_2$ be the TWO's response in $\G1(\obs,\Gbs)$ and so on.

This defines a strategy $\sigma$ for ONE in $\G1(\Omega_{\underline{0}},\Sigma_{\underline{0}})$. Since $\sigma$ is not a winning strategy for ONE, consider a $\sigma$-play $$\sigma(\emptyset),f_{B_1,U_1}, \sigma(f_{B_1,U_1}),f_{B_2,U_2},\dots$$ which is lost by ONE. So $\{f_n:n\in \nb\}\in \Sigma_{\underline{0}}$, where $f_n=f_{B_n,U_n}$, $n\in \nb$. Consider $\{U_n:n\in \nb\}$. Since $f_n^{-1}(1,1)\subseteq U_n$ for each $n\in \nb$, by Lemma \ref{L8}, $\{U_n:n\in \nb\}$) is a $\gbs$-cover of $X$.

The corresponding $\psi$-play in $\G1(\obs,\Gbs)$ is \[\psi(\emptyset), U_1, \psi(U_1), U_2,\dotsc~.\] Since $\{U_n:n\in \nb\}\in \Gbs$, $\psi$ is not a winning strategy for ONE in
$\G1(\obs,\Gbs)$.
\ep

%The first result is in the context of the topology of strong uniform convergence on $\bk$ of a well known result in connection with the $\alpha_i$-properties in $C_p(X)$ (see \cite{salpha1},\cite{salpha2}).
Next we present some applications of well known $\alpha_i$-properties with respect to the topology of strong uniform convergence on $\bk$.
In the following results we make use of the fact that a sequence $\{f_n:n\in \nb\}$ in $(C(X),\tau^s_\bk)$ converges to $\underline{0}$ with respect to $\tau^s_\bk$ if and only if $\{|f_n|:n\in \nb\}$ converges to $\underline{0}$ with respect to $\tau^s_\bk$.

\begin{Th}
\label{Tf2}
Let $\bk$ be a bornology on $X$ with closed base and $X$ is $\bs$-Lindel\"{o}f. The following statements are equivalent.

\noindent$(1)$ $(C(X),\tau^s_\bk)$ satisfies $\alpha_2(\Omega_{\underline{0}}, \Sigma_{\underline{0}})$.

\noindent$(2)$ $(C(X),\tau^s_\bk)$ satisfies $\alpha_3(\Omega_{\underline{0}},\Sigma_{\underline{0}})$.

\noindent$(3)$ $(C(X),\tau^s_\bk)$ satisfies $\alpha_4(\Omega_{\underline{0}},\Sigma_{\underline{0}})$.

\noindent$(4)$ $X$ satisfies $\S1(\obs,\Gbs)$.
\end{Th}
\bp
We only prove $(3)\Rightarrow (4)$ and $(4)\Rightarrow (1)$.

$(3)\Rightarrow (4)$. Let $\{\uc_n:n\in \nb\}$ be a sequence of open $\bs$-covers of $X$.
Enumerate each $\uc_n$ bijectively as $\{U^n_k:k\in \nb\}$. For each $n$ the collection $\vc_n=\{U^1_{k_1}\cap U^2_{k_2}\cap\dots \cap U^n_{k_n}:1<k_1<k_2<\cdots <k_n\}$ is an open $\bs$-cover of $X$. Now define $A_n=\{f\in C(X):$ there is a $V\in \vc_n$ with $f(X\setminus V)=\{1\}\}$ for each $n$. By Lemma \ref{koc}, $A_n\in \Omega_{\underline{0}}$. Now apply $(3)$ to $\{A_n:n\in \nb\}$ to obtain a sequence $n_1<n_2<\cdots$ and a $A\in \Sigma_{\underline{0}}$ such that $A_{n_i}\cap A\neq \emptyset$ for each $i\in \nb$. Let $f_{n_i}\in A_{n_i}\cap A$. Clearly $\{f_{n_i}:i\in \nb\}\in \Sigma_{\underline{0}}$. Now there is a $V^{n_i}_{l_i}\in \vc_{n_i}$ such that $f_{n_i}(X\setminus V^{n_i}_{l_i})=\{1\}$ for each $i$. We claim that $\{V^{n_i}_{l_i}:i\in \nb\}$ is a $\gbs$-cover of $X$. Let $B\in \bk$. Since $\{f_{n_i}:i\in \nb\}\in \Sigma_{\underline{0}}$, for the neighbourhood $[B,1]^s(\underline{0})$ of $\underline{0}$ there is a $i_0$ such that $f_{n_i}\in [B,1]^s(\underline{0})$ for all $i\geq i_0$ i.e. there is a $\delta_i>0$ such that $|f_{n_i}(x)|<1$ for all $x\in B^{\delta_i}$ and $i\geq i_0$ i.e. $B^{\delta_i}\subseteq f_{n_i}^{-1}(-1,1)\subseteq V^{n_i}_{l_i}$ for all $i\geq i_0$. This shows that $\{V^{n_i}_{l_i}:i\in \nb\}$ is a $\gbs$-cover of $X$. Now for each $n$ with $1\leq n \leq n_1$ let $U_n\in \uc_n$ be the $n$-th component in the representation of $V^{n_1}_{l_1}$. Also for each $i>1$ and each $n\in (n_{i-1},n_{i}]$ let $U_n\in \uc_n$ be the $n$-th component in the representation of $V^{n_{i}}_{l_{i}}$. Then it is easy to verify that $\{U_n:n\in \nb\}$ is a $\gbs$-cover of $X$. Hence $X$ satisfies $\S1(\obs,\Gbs)$.

$(4)\Rightarrow (1)$. Let $\{A_n:n\in \nb\}$ be a sequence of elements in $\Omega_{\underline{0}}$. Since $X$ is $\bs$-Lindel\"{o}f, $(C(X),\tau^s_\bk)$ has countable tightness. Therefore we can assume that $A_n$'s are countable. Say $A_n=\{f_{n,k}:k\in \nb\}$ for each $n$. We can also assume that $f_{n,k}\geq 0$ for each $n$ and $k$. By \cite[Theorem 3.7]{dcpdsd} and Theorem \ref{Tf1}, we can say that ONE has no winning strategy in $\G1(\Omega_{\underline{0}},\Sigma_{\underline{0}})$ on $(C(X),\tau^s_\bk)$. We now define a strategy $\sigma$ for ONE in $\G1(\Omega_{\underline{0}},\Sigma_{\underline{0}})$ as follows.

Define $\sigma(\emptyset)=A_1$, the first move of ONE. TWO responds with $f_{1,k_{i_1}}\in A_1$. Let $D_1=\{f_{1,k_1}+f_{2,k_2}:k_{i_1}<k_1<k_2\}$. It is easy to see that $D_1\in \Omega_{\underline{0}}$. Define $\sigma(f_{1,k_{i_1}})=D_1$. TWO responds with $f_{1,k_{i_2}}+f_{2,k_{i_3}}$. Again consider $D_2=\{f_{1,k_1}+f_{2,k_2}+f_{3,k_3}:k_{i_3}<k_1<k_2<k_3\}\in \Omega_{\underline{0}}$ and define $\sigma(f_{1,k_{i_1}}, f_{1,k_{i_2}}+f_{2,k_{i_3}})=D_2$. TWO responds with
$f_{1,k_{i_4}}+f_{2,k_{i_5}}+f_{3,k_{i_6}}$ and so on.

This defines a strategy $\sigma$ for ONE in $\G1(\Omega_{\underline{0}},\Sigma_{\underline{0}})$. Since $\sigma$ is not a winning strategy, consider a play
$$\sigma(\emptyset), f_{1,k_{i_1}}, \sigma(f_{1,k_{i_1}}), f_{1,k_{i_2}}+f_{2,k_{i_3}}, \sigma(f_{1,k_{i_1}}, f_{1,k_{i_2}}+f_{2,k_{i_3}}), f_{1,k_{i_4}}+f_{2,k_{i_5}}+f_{3,k_{i_6}}, \dotsc$$
which is lost by ONE. Therefore

$$f_{1,k_{i_1}}, f_{1,k_{i_2}}+f_{2,k_{i_3}}, f_{1,k_{i_4}}+f_{2,k_{i_5}}+f_{3,k_{i_6}},\dotsc$$ converges to $\underline{0}$ with respect to $\tau^s_\bk$. Since $f_{n,m}\geq 0$ for $n,m\in \nb$, the sequence

$$f_{1,k_{i_1}}, f_{1,k_{i_2}}, f_{2,k_{i_3}}, f_{1,k_{i_4}}, f_{2,k_{i_5}}, f_{3,k_{i_6}},\dotsc$$ also converges to $\underline{0}$ with respect to $\tau^s_\bk$, which contains infinitely many elements from each $A_n$. Hence $(1)$ holds.

\ep

Summarizing Theorem \ref{Tp3}, Theorem \ref{Tf1} and Theorem \ref{Tf2} we obtain the following.
\begin{Cor}
\label{C1}
Let $\bk$ be a bornology on $X$ with closed base and $X$ is $\bs$-Lindel\"{o}f. The following statements are equivalent.

\noindent$(1)$ ONE has no winning strategy in $\G1(\Omega_{\underline{0}},\Sigma_{\underline{0}})$ on $(C(X),\tau^s_\bk)$.

\noindent$(2)$ $(C(X),\tau^s_\bk)$ satisfies $\alpha_2(\Omega_{\underline{0}}, \Sigma_{\underline{0}})$.

\noindent$(3)$ $(C(X),\tau^s_\bk)$ satisfies $\alpha_3(\Omega_{\underline{0}},\Sigma_{\underline{0}})$.

\noindent$(4)$ $(C(X),\tau^s_\bk)$ satisfies $\alpha_4(\Omega_{\underline{0}},\Sigma_{\underline{0}})$.

\noindent$(5)$ $X$ satisfies $\S1(\obs,\Gbs)$.

\noindent$(6)$ ONE has no winning strategy in $\G1(\obs,\Gbs)$ on $X$.

\noindent$(7)$ $X$ satisfies $\obs\rightarrow (\Gbs)^n_k$ for each $n,k\in \nb$, provided $X$ is $\bs$-Lindel\"{o}f.

\end{Cor}

\begin{Rem}
Moreover using \cite[Corollary 2.10]{cmk}, the following equivalent conditions can also be added to Corollary \ref{C1}.

\noindent$(8)$ $(C(X),\tau^s_\bk)$ is Fr\'{e}chet-Urysohn.

\noindent$(9)$ $(C(X),\tau^s_\bk)$ is strictly Fr\'{e}chet-Urysohn.

\noindent$(10)$ Every open $\bs$-cover of $X$ contains a countable set which is a $\gbs$-cover of $X$.
\end{Rem}

\begin{Th}
\label{Tf}
Let $\bk$ be a bornology on $X$ with closed base. The following statements are equivalent.\\
\noindent$(1)$ $(C(X),\tau^s_\bk)$ satisfies $\alpha_2(\Sigma_{\underline{0}}, \Sigma_{\underline{0}})$.\\
\noindent$(2)$ $(C(X),\tau^s_\bk)$ satisfies $\alpha_3(\Sigma_{\underline{0}}, \Sigma_{\underline{0}})$.\\
\noindent$(3)$ $(C(X),\tau^s_\bk)$ satisfies $\alpha_4(\Sigma_{\underline{0}}, \Sigma_{\underline{0}})$.\\
\noindent$(4)$ $(C(X),\tau^s_\bk)$ satisfies $\S1(\Sigma_{\underline{0}}, \Sigma_{\underline{0}})$.\\
\noindent$(5)$ ONE has no winning strategy in $\G1(\Sigma_{\underline{0}}, \Sigma_{\underline{0}})$.
\end{Th}
\bp

We only give  proof of the implications $(3)\Rightarrow (4)$, $(4)\Rightarrow (5)$ and $(4)\Rightarrow (1)$.

$(3)\Rightarrow (4)$. Let $\{S_n:n\in \nb\}$ be a sequence of elements in $\Sigma_{\underline{0}}$ where $S_n=\{f_{n,m}:m\in \nb\}$. We assume that $f_{n,m}\geq 0$ for all $n,m\in \nb$. Fix a $n\in \nb$, we construct a new sequence $g_{n,m}=f_{1,m}+f_{2,m}+ \cdots +f_{n,m}$, $m\in \nb$. Then the sequence $\{g_{n,m}:m\in \nb\}$ converges to $\underline{0}$ with respect to $\tau^s_\bk$. Applying $(3)$ to the sequence $\{T_n:n\in \nb\}$, where $T_n=\{g_{n,m}:m\in \nb\}$ for each $n\in \nb$, we obtain an increasing sequence  $1=n_0<n_1<n_2<\cdots$ of positive integers such that $\{g_{{n_i},{m_i}}:i\in \nb\}$ converges to $\underline{0}$ and $g_{{n_i},{m_i}}\in T_{n_i}$ for each $i$. For each integer $k\geq 0$ and $j\in \nb$ with $n_k<j\leq  n_{k+1}$, note that $g_{n_{k+1},m_{k+1}}=f_{1,m_{k+1}}+\cdots +f_{n_{k+1},m_{k+1}}$, and now choose $h_j= f_{j,m_{k+1}}$. We will show that $\{h_j:j\in \nb\}$ where $h_j\in S_j$ witnesses $\S1(\Sigma_{\underline{0}}, \Sigma_{\underline{0}})$. For this we need to show that the sequence $\{h_j:j\in \nb\}$ converges to $\underline{0}$ with respect to $\tau^s_\bk$.

Let $[B,\varepsilon]^s(\underline{0})$ ($B\in \bk, \varepsilon>0$) be a neighbourhood of $\underline{0}$. There exists a $i_0\in \nb$ such that $g_{n_{i+1},m_{i+1}}\in [B,\varepsilon]^s(\underline{0})$ for all $i\geq  i_0$. This means that there is a $\delta_i>0$ such that $|f_{1,m_{i+1}}(x)+\cdots +f_{n_{i+1},m_{i+1}}(x)|<\varepsilon$ for all $x\in B^{\delta_i}$ and for all $i\geq i_0$ i.e. $|f_{j,m_{i+1}}(x)|<\varepsilon$ for all $x\in B^{\delta_i}$ and for all $j=1,\dots,n_{i+1}$, $i\geq i_0$ i.e. $f_{j,m_{i+1}}\in [B,\varepsilon]^s(\underline{0})$ for all $n_i<j\leq n_{i+1}$ and for all $i\geq i_0$ i.e. $h_j\in [B,\varepsilon]^s(\underline{0})$ for all $j>n_{i_0}$. Therefore $\{h_j:j\in \nb\}$ converges to $\underline{0}$ with respect to $\tau^s_\bk$. Hence $(4)$ holds.\\

$(4)\Rightarrow (5)$. Let $\sigma$ be a strategy for ONE in $\G1(\Sigma_{\underline{0}}, \Sigma_{\underline{0}})$. The first move of ONE is $\sigma(\emptyset)$. Enumerating $\sigma(\emptyset)$ bijectively as $\{f_{(n)}:n\in \nb\}$. Assume that for each finite sequence of natural numbers $\tau$ of length at most $m$, $U_\tau$ is already defined. Now define $\{f_{(n_1,\dots,n_m,k)}:k\in \nb\}$ to be $\sigma(f_{(n_1)},\dots, f_{(n_1,\dots,n_m)})\setminus \{f_{(n_1)},\dots, f_{(n_1,\dots,n_m)}\}$. We have for each finite sequence of natural numbers $\tau$, $\{f_{\tau\frown n}:n\in \nb\}\in \Sigma_{\underline{0}}$. Apply $\S1(\Sigma_{\underline{0}}, \Sigma_{\underline{0}})$ to choose a $n_\tau$ such that $\{U_{\tau\frown n_\tau}:\tau$ is a finite sequence of natural numbers$\}$ converges to $\underline{0}$ with respect to $\tau^s_{\bk}$. Inductively define a sequence $n_1=n_{\emptyset}$, $n_{k+1}=n_{(n_1,\dots, n_k)}$ for $k\geq 1$. Then the sequence
\[U_{n_1}, U_{(n_1,n_2)},\dotsc, U_{(n_1,\dots,n_k)},\dotsc\]
converges to $\underline{0}$ with respect to $\tau^s_{\bk}$. Since it is actually a sequence of moves of TWO during a play in $\G1(\Sigma_{\underline{0}}, \Sigma_{\underline{0}})$, $\sigma$ is not a winning strategy for ONE.

$(4)\Rightarrow (1)$. Let $\{S_n:n\in \nb\}$ be a sequence of elements in $\Sigma_{\underline{0}}$. From each $S_n$ we can easily make new sequences $S_{n,m}$ which are pairwise disjoint and each of which converges to $\underline{0}$. Consider the collection $\{S_{n,m}:n,m\in \nb\}$. Now apply $\S1(\Sigma_{\underline{0}}, \Sigma_{\underline{0}})$ to  choose a $f_{n,m}\in S_{n,m}$ for $n,m\in \nb$ such that $\{f_{n,m}:n,m\in \nb\}$ converges to $\underline{0}$. Clearly the sequence $\{f_{n,m}:n,m\in \nb\}$ contains infinitely many elements from each $S_n$. Hence $(1)$ holds.
\ep

%\begin{Th}
%\label{Tgn2}
%Let $\bk$ be bornology on $X$ with closed base. Then the following statements are equivalent.\\
%\noindent$(1)$ $X$ has the $\bs$-Gerlits-Nagy property.\\
%\noindent$(2)$ $(C(X),\tau_\bk^s)$ satisfies $\S1(\Omega_{\underline{0}},\Omega_{\underline{0}}^{gp})$.
%\end{Th}
%\bp
%$(1)\Rightarrow (2)$. $X$ satisfies $\S1(\obs,\obs^{gp})$ implies that $X$ satisfies $\S1(\obs,\obs)$ which in turn implies that  $(C(X),\tau_\bk^s)$ satisfies $\S1(\Omega_{\underline{0}},\Omega_{\underline{0}})$ by \cite[Theorem 2.3]{cmk}. Now it remains to show that every countable element in $\Omega_{\underline{0}}$ is groupable. Clearly $X$ is $\bs$-Lindel\"{o}f and every countable open $\bs$-cover of $X$ is $\bs$-groupable. By \cite[Theorem 5.3]{dcpdsd}, every element in $\Omega_{\underline{0}}$ is groupable. Hence $(C(X),\tau_\bk^s)$ satisfies $\S1(\Omega_{\underline{0}},\Omega_{\underline{0}}^{gp})$.\\

%$(2)\Rightarrow (1)$. $(C(X),\tau_\bk^s)$ satisfies $\S1(\Omega_{\underline{0}},\Omega_{\underline{0}}^{gp})$ implies that $X$ satisfies $\S1(\obs,\obs)$ by \cite[Theorem 2.3]{cmk}. Using the fact that every elements in $\Omega_{\underline{0}}$ is groupable and applying \cite[Theorem 5.3]{dcpdsd}, we conclude that every countable open $\bs$-cover of $X$ is $\bs$-groupable. So $X$ satisfies $\S1(\obs,\obs^{gp})$. Hence $(1)$ holds.
%\ep

\subsection{Some observations on $X^n$ and $(C(X),\tau^s_\bk)$}
\begin{Th}
\label{Tfp-1}
Let $\bk$ be a bornology on $X$ with compact base. The following statements are equivalent.\\
\noindent$(1)$ $(C(X),\tau^s_{\bk})$  has countable tightness.\\
\noindent$(2)$ $X^n$ is $\bns$-Lindel\"{o}f for each $n\in \nb$.
\end{Th}
\bp
$(1)\Rightarrow (2)$. Let $\uc$ be an open $\bns$-cover of $X^n$. By Lemma \ref{Lpdt}, there is an open $\bs$-cover $\vc$ of $X$ such that $\{V^n:V\in \vc\}$ refines $\uc$. Now for $B\in \bk$ there exist a $\delta>0$ and a $V\in \vc$ such that $B^{2\delta}\subseteq V$. Let $\vc_B=\{V\in \vc:B^{2\delta}\subseteq V\}$. For each $V\in \vc_B$ choose a $f_{B,V}\in C(X)$ such that $f_{B,V}(B^\delta)=\{0\}$ and $f_{B,V}(X\setminus V)=\{1\}$. Consider the set $F=\{f_{B,V}:B\in \bk, V\in \vc_B\}$. Clearly $\underline{0}\in \overline{F}$. By $(1)$, there is a countable subset $F'$ of $F$ such that $\underline{0}\in \overline{F'}$. Let $\wc=\{V:f_{B,V}\in F'\}$. Then $\wc$ is a countable subset of $\vc$. Now for $V\in \wc$ choose a $U\in \uc$ such that $V^n\subseteq U$. Consider the set $\zc=\{U:V^n\subseteq U$ for $V\in \wc\}$. Clearly $\zc$ is a countable subset of $\uc$. We show that $\zc$ is an open $\bns$-cover of $X^n$. Let $B^n\in \bk^n$. As $[B,1]^s(\underline{0})\cap F'\neq \emptyset$, there is a $f_{B_1,V_1}\in F'$ with $f_{B_1,V_1}\in [B,1]^s(\underline{0})$ i.e. there exists a $\delta>0$ such that $B^\delta\subseteq f^{-1}_{B_1,V_1}(-1,1)\subseteq V_1$ i.e. $(B^\delta)^n\subseteq V^n_1\subseteq U_1$ for some $U_1\in \zc$ i.e. $(B^n)^\delta\subseteq U_1$. So $\zc$ is a countable $(\bk^n)^s$-subcover of $\uc$.

$(2)\Rightarrow (1)$. Let $A\subseteq (C(X),\tau^s_{\bk^s})$ with $\underline{0}\in \overline{A}$. Let $\uc_m=\{g^{-1}(-\frac{1}{m},\frac{1}{m}):g\in A\}$ for $m\in \nb$. Then by Lemma \ref{koc} $\uc_m$ is an open $\bs$-cover of $X$ for each $m\in \nb$. Let $\vc_m=\{U^n:U\in \uc_m\}$. Clearly $\vc_m$ is an open $\bns$-cover of $X^n$. By $(2)$, for each $m\in \nb$ there is a countable $\bns$-subcover $\wc_m=\{U^n_{k,m}:k\in \nb\}$ of $\vc_m$, where $U_{k,m}=g^{-1}_{k,m}(-\frac{1}{m},\frac{1}{m})$ for $k\in \nb$. Choose $A'=\{g_{k,m}:k,m\in \nb\}$, a countable subset of $A$. We show that $\underline{0}\in \overline{A'}$. Let $[B,\varepsilon]^s(\underline{0})$ be a neighbourhood of $\underline{0}$, where $B\in \bk$ and $\varepsilon>0$. Choose a $m\in \nb$ with $\frac{1}{m}<\varepsilon$. Now for $B^n\in \bn$ there exist a $\delta>0$ and a $U^n_{k,m}\in \wc_m$ such that $(B^n)^\delta\subseteq U^n_{k,m}$ i.e. $B^\delta\subseteq U_{k,m}= g^{-1}_{k,m}(-\frac{1}{m},\frac{1}{m})$ i.e. $g_{k,m}\in [B,\varepsilon]^s(\underline{0})$ i.e. $[B,\varepsilon]^s(\underline{0})\cap A'\neq \emptyset$. So $\underline{0}\in \overline{A'}$. Hence $(1)$ holds.
\ep

 The next three results can be obtained by applying Theorem \ref{Tpdt1} to each \cite[Theorem 2.3]{cmk}, \cite[Theorem 2.5]{cmk} and \cite[Theorem 2.7]{cmk}  respectively.
\begin{Th}
\label{Tfp-2}
Let $\bk$ be a bornology on $X$ with compact base. The  following statements are equivalent.

$(1)$ $(C(X),\tau^s_\bk)$ has countable strong fan tightness.

$(2)$ $X^n$ satisfies $\S1(\obns,\obns)$ for each $n\in \nb$.
\end{Th}

\begin{Th}
\label{Tfp-3}
Let $\bk$ be a bornology on $X$ with compact base. The  following statements are equivalent.

$(1)$ $(C(X),\tau^s_\bk)$ has countable fan tightness.

$(2)$ $X^n$ satisfies $\Sf(\obns,\obns)$ for each $n\in \nb$.
\end{Th}

\begin{Th}
\label{Tfp-4}
Let $\bk$ be a bornology on $X$ with compact base. The  following statements are equivalent.

$(1)$ $(C(X),\tau^s_\bk)$ is strictly Fr\'{e}chet-Urysohn.

$(2)$ $X^n$ satisfies $\S1(\obns,\Gamma_{\bns})$ for each $n\in \nb$.
\end{Th}

In line of \cite[Theorem 21]{coocvii}, we can obtain a bornological version which presents similar characterization using the topology of strong uniform convergence $\tau^s_\bk$ on $\bk$.
\begin{Th}
\label{Tfp-5}
Let $\bk$ be a bornology on $X$ with closed base. The  following statements are equivalent.

$(1)$ $X$ has the $\bs$-Hurewicz property.

$(2)$ $(C(X),\tau_\bk^s)$ has countable fan tightness and Reznichenko's property i.e. it satisfies $\Sf(\Omega_{\underline{0}},\Omega^{gp}_{\underline{0}})$.

$(3)$ ONE does not have a winning strategy in $\Gf(\Omega_{\underline{0}},\Omega^{gp}_{\underline{0}})$ on $(C(X),\tau_\bk^s)$.
\end{Th}
\bp
The equivalence $(1)\Leftrightarrow (2)$ is due to \cite[Theorem 5.8]{dcpdsd}. The implication $(3)\Rightarrow (2)$ is easily followed. We prove the implication $(1)\Rightarrow (3)$.

$(1)\Rightarrow (3)$. $X$ has the $\bs$-Hurewicz property implies that ONE does not have a winning strategy in the $\bs$-Hurewicz game on $X$. Let $\sigma$ be a strategy for ONE in the game $\Gf(\Omega_{\underline{0}},\Omega^{gp}_{\underline{0}})$. We use $\sigma$ to define a strategy $\psi$ for ONE in the $\bs$-Hurewicz game as follows.

The first move of ONE in $\Gf(\Omega_{\underline{0}},\Omega^{gp}_{\underline{0}})$ is $\sigma(\emptyset)$. Let $\uc_1=\{f^{-1}(-1,1):f\in \sigma(\emptyset)\}$ and assume that $X\not\in \uc_1$. By Lemma \ref{koc}, $\uc_1$ is an open $\bs$-cover of $X$. Define $\psi(\emptyset)=\uc_1$, the first move of ONE the $\bs$-Hurewicz game. TWO responds by choosing a finite subset $\vc_1$ of $\psi(\emptyset)$. Let $C_1\subseteq \sigma(\emptyset)$ be a finite collection of functions such that $\vc_1=\{f^{-1}(-1,1):f\in C_1\}$. Then $C_1$ is a legitimate move of TWO in $\Gf(\Omega_{\underline{0}},\Omega^{gp}_{\underline{0}})$.

To define $\psi(\vc_1)$, look at the move $\sigma(C_1)$ of ONE in $\Gf(\Omega_{\underline{0}},\Omega^{gp}_{\underline{0}})$. Let $A_2=\sigma(C_1)\setminus C_1$. Clearly $A_2\in \Omega_{\underline{0}}$ and $\uc_2=\{f^{-1}(-\frac{1}{2},\frac{1}{2}):f\in A_2\}$ is an open $\bs$-cover of $X$ again by Lemma \ref{koc}. Now define $\psi(\vc_1)=\uc_2$. TWO responds by choosing a finite subset $\vc_2$ of $\psi(\vc_1)$. Let $C_2$ be the finite subset of $A_2$ with $\vc_2=\{f^{-1}(-\frac{1}{2},\frac{1}{2}):f\in C_2\}$. Thus $C_2$ is a legitimate move of TWO in $\Gf(\Omega_{\underline{0}},\Omega^{gp}_{\underline{0}})$. Again look at $\sigma(C_1,C_2)$ and let $A_3=\sigma(C_1,C_2)\setminus \{C_1,C_2\}$ which is in $\Omega_{\underline{0}}$. $\uc_3=\{f^{-1}(-\frac{1}{3},\frac{1}{3}):f\in A_3\}$ is an open $\bs$-cover of $X$. Define $\psi(\vc_1,\vc_2)=\uc_3$. TWO responds by choosing a finite set $\vc_3\subseteq \psi(\vc_1,\vc_2)$. Let $C_3\subseteq A_3$ be a finite set with $\vc_3=\{f^{-1}(-\frac{1}{3},\frac{1}{3}):f\in C_3\}$ and so on.

This defines a strategy $\psi$ for ONE in the $\bs$-Hurewicz game. Since $\psi$ is not a winning strategy, consider a $\psi$-play
$$\psi(\emptyset), \vc_1, \psi(\vc_1), \vc_2, \psi(\vc_1,\vc_2), \dots$$ which is lost by ONE. Therefore for $B\in \bk$ there exist a $n_0\in \nb$ and a sequence $\{\delta_n:n\geq n_0\}$ of positive real numbers satisfying $B^{\delta_n}\subseteq U$ for some $U\in \vc_n$ for all $n\geq n_0$.

Now corresponding to the $\psi$-play there is a $\sigma$-play in $\Gf(\Omega_{\underline{0}},\Omega^{gp}_{\underline{0}})$
$$\sigma(\emptyset), C_1, \sigma(C_1), C_2, \sigma(C_1,C_2), C_3, \dots$$
We will show that $\cup_{n\in \nb}C_n\in \Omega^{gp}_{\underline{0}}$.

It is clear from the construction of the game that $C_n$'s are pairwise disjoint and for $U\in \vc_n$, $U=f^{-1}(-\frac{1}{n},\frac{1}{n})$ for some $f\in C_n$. Let $[B,\varepsilon]^s(\underline{0})$ be a neighbourhood of $\underline{0}$ where $B\in \bk$ and $\varepsilon>0$. For $B\in \bk$ there is a $n_0\in \nb$ and a sequence $\{\delta_n:n\geq n_0\}$ of positive real numbers satisfying $B^{\delta_n}\subseteq U$ for some $U\in \vc_n$ for all $n\geq n_0$. Choose a $n_1\in \nb$ such that $\frac{1}{n_1}<\varepsilon$. Choose $n_2=\max\{n_0,n_1\}$. Then $B^{\delta_n}\subseteq f^{-1}(-\frac{1}{n},\frac{1}{n})\subseteq f^{-1}(-\varepsilon,\varepsilon)$ for some $f\in C_n$ for all $n\geq n_2$ i.e. $f\in [B,\varepsilon]^s(\underline{0})$ for some $f\in C_n$ for all $n\geq n_2$ i.e. $[B,\varepsilon]^s(\underline{0})\cap C_n\neq \emptyset$ for all $n\geq n_2$. Therefore $\{C_n:n\in \nb\}$ witnesses the groupability of $\cup_{n\in \nb}C_n$. Hence $\sigma$ is not a winning strategy for ONE in $\Gf(\Omega_{\underline{0}},\Omega^{gp}_{\underline{0}})$.
\ep

Theorem \ref{Tfp-5} together with Theorem \ref{Tpdt2} gives the following.
\begin{Cor}
\label{C2}
Let $\bk$ be a bornology on $X$ with compact base. The  following statements are equivalent.

$(1)$ $X^n$ has the $\bns$-Hurewicz property for each $n\in \nb$.

$(2)$ $(C(X),\tau_\bk^s)$ has countable fan tightness and Reznichenko's property i.e. it satisfies $\Sf(\Omega_{\underline{0}},\Omega^{gp}_{\underline{0}})$.

$(3)$ ONE does not have a winning strategy in $\Gf(\Omega_{\underline{0}},\Omega^{gp}_{\underline{0}})$ on $(C(X),\tau_\bk^s)$.
\end{Cor}

\begin{Th}
\label{Tfp-6}
Let $\bk$ be a bornology on $X$ with closed base. The  following statements are equivalent.

$(1)$ $X$ has the $\bs$-Gerlits-Nagy property.

$(2)$ $(C(X),\tau_\bk^s)$ has countable strong fan tightness and Reznichenko's property.
\end{Th}
\bp
$(1)\Rightarrow (2)$. By Theorem \ref{Tgn-1}, $X$ has the $\bs$-Hurewicz property as well as it satisfies $\S1(\obs,\obs)$. Also from \cite[Theorem 5.8]{dcpdsd}, $(C(X),\tau_\bk^s)$ has the Reznichenko's property.  Again $X$ satisfies $\S1(\obs,\obs)$ implies that $(C(X),\tau_\bk^s)$ has countable strong fan tightness by \cite[Theorem 2.3]{cmk}. Hence $(2)$ holds.

$(2)\Rightarrow (1)$. Using \cite[Theorem 5.8]{dcpdsd} and  \cite[Theorem 2.3]{cmk}, we have $X$ has the $\bs$-Hurewicz property as well as it satisfies $\S1(\obs,\obs)$. Hence $X$
the $\bs$-Gerlits-Nagy property.
\ep

\begin{Cor}
\label{C3}
Let $\bk$ be a bornology on $X$ with compact base. The  following statements are equivalent.

$(1)$ $X^n$ has the $\bns$-Gerlits-Nagy property for each $n\in \nb$.

$(2)$ $(C(X),\tau_\bk^s)$ has countable strong fan tightness and Reznichenko's property.
\end{Cor}

\begin{Th}
\label{Th&gn}
Let $\bk_1, \bk_2$ be bornologies on $X_1$ and $X_2$ respectively with closed bases. If $(C(X_1),\tau^s_{\bk_1})$ and $(C(X_2),\tau^s_{\bk_2})$ are homeomorphic, then the following statements are true.

$(1)$ $X_1$ has the $\bk^s_1$-Hurewicz property if and only if $X_2$ has the $\bk^s_2$-Hurewicz property.

$(2)$ $X_1$ has the $\bk^s_1$-Gerlits-Nagy property if and only if $X_2$ has the $\bk^s_2$-Gerlits-Nagy property.
\end{Th}
\bp
We only prove $(1)$ as the proof of $(2)$ is analogous.

$(1)$. Let $\psi:(C(X_1),\tau^s_{\bk_1})\rightarrow (C(X_2),\tau^s_{\bk_2})$ be a homeomorphism. Suppose that $X_1$ has the $\bk^s_1$-Hurewicz property. By \cite[Theorem 5.8]{dcpdsd}, $(C(X_1),\tau^s_{\bk_1})$ satisfies $\Sf(\Omega_{\underline{0}_1},\Omega_{\underline{0}_1}^{gp})$. We will show that $(C(X_2),\tau^s_{\bk_2})$ satisfies  $\Sf(\Omega_{\underline{0}_2},\Omega_{\underline{0}_2}^{gp})$. For this we first show that $\Sf(\Omega_{\underline{0}_2},\Omega_{\underline{0}_2})$ holds and then prove that every countable element in $\Omega_{\underline{0}_2}$ is groupable. ($\underline{0}_1$, $\underline{0}_2$ are zero elements in $C(X_1)$ and $C(X_2)$ respectively).

Let $\{A_n:n\in \nb\}$ be a sequence of elements in $\Omega_{\underline{0}_2}$. It is easy to observe that if $A_n\in \Omega_{\underline{0}_2}$ then $\psi^{-1}(A_n)\in \Omega_{\underline{0}_1}$. Apply $\Sf(\Omega_{\underline{0}_1},\Omega_{\underline{0}_1})$ to $\{\psi^{-1}(A_n):n\in \nb\}$ to choose a finite subset $B_n$ of $\psi^{-1}(A_n)$ for each $n\in \nb$ such that $\cup_{n\in \nb} B_n\in \Omega_{\underline{0}_1}$. Now $\psi(B_n)$ is a finite subset of $A_n$ for each $n\in \nb$. Take a neighbourhood $[B,\varepsilon]^s(\underline{0}_2)$ of $\underline{0}_2$ where $B\in \bk_2$, $\varepsilon>0$. Consequently there is a neighbourhood $[B',\delta]^s(\underline{0}_1)$ of $\underline{0}_1$ where $B'\in \bk_1$, $\delta>0$ such that $\psi([B',\delta]^s{(\underline{0}_1)})\subseteq [B,\varepsilon]^s(\underline{0}_2)$. Again $[B',\delta]^s(\underline{0}_1)\cap (\cup_{n\in \nb} B_n)\neq \emptyset$ i.e. $[B,\varepsilon]^s(\underline{0}_2)\cap (\cup_{n\in \nb} \psi(B_n))\neq \emptyset$. So $\cup_{n\in \nb}\psi(B_n)\in \Omega_{\underline{0}_2}$. Hence $\Sf(\Omega_{\underline{0}_2},\Omega_{\underline{0}_2})$ holds.

Finally let $A$ be a countable element in $ \Omega_{\underline{0}_2}$. Clearly $\psi^{-1}(A)\in \Omega_{\underline{0}_1}$. Since every countable element in $\Omega_{\underline{0}_1}$ is groupable, there is sequence $\{C_n:n\in \nb\}$ of pairwise disjoint finite subsets of $\psi^{-1}(A)$ witnessing the groupability of $\psi^{-1}(A)$. It is easy to verify that $\{\psi(C_n):n\in \nb\}$ is a sequence of pairwise disjoint finite subsets of $A$ witnessing the groupability of $A$. Hence  $(C(X_2),\tau^s_{\bk_2})$ satisfies  $\Sf(\Omega_{\underline{0}_2},\Omega_{\underline{0}_2}^{gp})$ and  $Y$ has the $\bk^s_2$-Hurewicz property by \cite[Theorem 5.8]{dcpdsd}.

Proof of the converse part follows with similar arguments.
\ep

{}
\end{document}